\documentclass[12pt]{article}
\usepackage{amsmath,indentfirst,enumitem,amsthm,amssymb,tikz,color,geometry,float,ragged2e,mathtools,graphicx,url,bigstrut,bm,booktabs,caption,subcaption}
\usepackage[utf8]{inputenc}
\usepackage{color}
\usepackage{amsfonts}
\usepackage{mathtools}
\usetikzlibrary{positioning,calc,arrows,automata}
\newtheorem{proposition}{PROPOSITION}

\theoremstyle{definition}
\newtheorem{definition}{DEFINITION}
\newtheorem{remark}{REMARK}
\newtheorem{example}{EXAMPLE}[section]
\theoremstyle{plain}
\newcommand{\R}{\mathbb{R}}
\newcommand{\C}{\mathbb{C}}

\newcommand{\nc}{\color{black}}

\author{Silvia Noschese\thanks{Dipartimento di Matematica, SAPIENZA Universit\`a di Roma,
P.le Aldo Moro 5, 00185 Roma, Italy. Email: \texttt{noschese@mat.uniroma1.it} }
\and
Lothar Reichel\thanks{Department of Mathematical Sciences, Kent State University, Kent,
OH 44242, USA. Email: \texttt{reichel@math.kent.edu}.}
}

\title{Sensitivity of Perron and Fiedler eigenpairs to structural perturbations of a
network}
\date{}

\begin{document}
\setlength{\arrayrulewidth}{1pt}
\maketitle

\begin{abstract}
One can estimate the change of the Perron and Fiedler values for a connected network when
the weight of an edge is perturbed by analyzing relevant entries of the Perron and Fiedler
vectors. This is helpful for identifying edges whose weight perturbation causes the 
largest change in the Perron and Fiedler values. It also is important to investigate the 
sensitivity of the Perron and Fiedler vectors to perturbations. Applications of the
perturbation analysis include the identification of edges that are critical for the 
structural robustness of the network.
\end{abstract}

\providecommand{\keywords}[1]
{\small\textbf{Keywords:} #1 }

\section{Introduction}\label{sec1}
Complex systems of interacting components often can be modeled by a graph 
${\cal G}=\{N,E\}$ that is made up of a set of $n$ nodes $N=\{1,2,\dots,n\}$ and a set of 
$m$ edges $E$. The nodes are connected by edges and, typically, $1\ll m\ll n^2$. Examples 
of graphs include social networks and transportation networks, such as road and airline 
networks; see, e.g., Estrada \cite{Es} and Newman \cite{Ne} for many examples of networks.
We will refer to ${\cal G}$ both as a graph and as a network.

This paper considers connected simple undirected weighted graphs, i.e., graphs that do not
have self-loops or multiple edges between adjacent nodes, and all edge weights are 
positive. The \emph{degree} of a node is the sum of the weights of the edges that start at 
the node. We say that node $i$ is connected to node $j$ if there is an edge between these
nodes. Let $A=[a_{ij}]_{i,j=1}^n\in\R^{n\times n}$ denote the \emph{adjacency matrix}
associated with ${\cal G}$. Then $a_{ij}>0$ if node $i$ is connected to node $j$; the
entry $a_{ij}$ is referred to as the \emph{weight} of the edge from node $i$ to node $j$.
If there is no edge between these nodes, then $a_{ij}=0$. The fact that ${\cal G}$ is 
undirected means that $a_{ij}=a_{ji}$ for all $1\leq i,j\leq n$. Thus, $A$ is symmetric.
The absence of self-loops implies that all diagonal entries of $A$ vanish. It follows from
the fact that the graph ${\cal G}$ is connected that the matrix $A$ is irreducible. A
graph is said to be unweighted if all nonvanishing entries of $A$ are one.

Let $\rho$ denote the spectral radius, also known as the \emph{Perron value}, of $A$ and 
let $u=[u_1,u_2,\ldots,u_n]^T\in \R^{n}$ be the unique associated eigenvector of unit 
Euclidean norm with all entries positive associated with $\rho$. Here and throughout this
paper the superscript $(\cdot)^T$ denotes transposition. The vector $u$ is known as the 
\emph{Perron vector} of $A$. The fact that $\rho$ is a simple eigenvalue and that $u$ can
be scaled to have positive entries only follows from the Perron-Frobenius Theorem; see, 
e.g., \cite[Section 2]{PF} for discussions. The Perron vector plays an important role in 
the ranking of nodes of a network; see \cite{Bo}.

The graph Laplacian associated with $\cal G$ is the symmetric matrix $L=D-A$, where 
$D={\rm diag}[d_1,d_2,\ldots,d_n]\in\R^{n\times n}$ is the degree matrix, i.e., 
$d_j=\sum_{k=1}^n a_{jk}$ is the degree of node $j$. It follows, e.g., by using 
Gershgorin disks, that the matrix $L$ is positive semidefinite. Let 
$1_n=[1,1,\ldots,1]^T\in\R^n$. Then $D1_n=A1_n$. Therefore, $L1_n=D1_n-A1_n=0$. Thus, 
$1_n$ is in the null space of $L$. Since ${\cal G}$ is connected, the null space of $L$ is
of dimension $1$; see \cite[Section 2]{PF}.

Let $\mu>0$ denote the second smallest eigenvalue of $L$ and let 
$v=[v_1,v_2,\ldots,v_n]^T\in\R^n$ be a corresponding eigenvector of unit Euclidean norm.
The eigenvalue $\mu$ is known as the \emph{Fiedler value} and the vector $v$ as a
\emph{Fiedler vector} of $L$. 
Eigenvectors of $L$ associated with distinct 
eigenvalues are orthogonal. Therefore, $v^T1_n=0$. It follows that the vector $v$ has both 
positive and negative entries.

The Fiedler value $\mu$ also is known as the \emph{algebraic connectivity} of the graph 
$\cal G$. It indicates how well-connected the graph $\cal G$ is; see Fiedler \cite{F1}. 
The signs of the entries of $v$ can be used to partition $\cal G$ into two connected subgraphs 
with node sets $N_1$ and $N_2$ such that $N=N_1\cup N_2$ and $N_1\cap N_2=\emptyset$. We
would like to determine a partition such that the node sets $N_1$ and $N_2$ have 
approximately the same number of nodes. An approach to seek to determine such a partition
of $\cal G$ is what we will refer to as the \emph{Fiedler procedure}. It reads: 
\begin{enumerate}
\item For $i=1,2,\ldots,n$: 
\item ~~~if $v_i<0$, then put node $i$ into $N_1$, else put $i$ into $N_2$;
\end{enumerate}
see \cite{F2}. 

The Fiedler procedure is quite effective  in minimizing the total {    cut weight, i.e., 
the sum of the weights of the edges to be removed in order to obtain a bipartition of the
graph}. 
We note that the squared component $v_i^2$ can 
be interpreted as the probability that node $i$ is in the component of the bipartition 
determined by the sign of $v_i$. When the bipartition is required to have subgraphs with 
the same number of nodes, it is the median $\beta$ of the entries of $v$ that plays the 
role of the watershed, i.e., we put node $i$ into the node set $N_1$ if $v_i<\beta$ 
and we put node $i$ into $N_2$ when $v_i\geq\beta$.

The accuracy, and therefore the usefulness, of the computed Perron and Fiedler values of
$A$, as well as of the associated eigenvectors, can be estimated from the sensitivity of the 
eigenpairs $(\rho,u)$ and $(\mu,v)$ to perturbations of the network that leave ${\cal G}$
connected and undirected. For instance, if the sensitivity of the Perron vector to small 
changes in the weights is negligible, then the removal of edges suggested by the Perron
vector (as described below) in a quest to simplify the network, also would be meaningful
if the weights would be slightly perturbed. This is important, because it may be difficult
to determine edge weights accurately.

A careful structured sensitivity analysis of the Perron and Fiedler values can be carried 
out by considering targeted perturbations of a single edge weight. Following the approach 
in \cite{MSN}, one can determine estimates of the relative variation of Perron and Fiedler 
values by using the entries of the Perron and Fiedler vectors relevant to the nodes that
are connected by the edge. We remark that the Fiedler value $\mu$ should not be viewed as 
a ``strict'' disconnectivity metric, because it may be easier to disconnect a graph with a
larger Fiedler value than a graph with a smaller one; see \cite[Section 4.31]{VM}.

The Perron value $\rho$ allows us to draw conclusions about the spread of a disease in a
network. Epidemiological theory predicts that viral contamination in a network dies out 
over time when the infection rate of a virus in an epidemic is smaller than the reciprocal
of the Perron value of the adjacency matrix for the network. In particular, to reduce the
spread of an infection in a network, one should decrease the weight of the edges whose 
weight reduction or removal leads to the largest decrease of the spectral radius; see, e.g., 
\cite{JKVV,MSN,NR1}. Removing an edge corresponds to setting the weight of this edge to 
zero.

It is the purpose of this paper to discuss how a network can be simplified by taking 
the Perron and Fielder values, as well as the conditioning of the Perron and Fiedler
vectors, into account. {\nc We remark that there are several other approaches available 
for measuring the significance of an edge and determining which edge(s) can be removed.
These approaches differ in how communication and the importance of edges in a graph is 
measured. For instance Altafini et al. \cite{ABCMP} measure the importance of edges with
the aid of the Kemeny constant, Arrigo and Benzi \cite{AB} consider edges between nodes
with (relatively) large subgraph centrality important, and Bergermann and Stoll \cite{BS}
are concerned with multiplex networks. This kind of networks also are discussed in 
\cite{NR3,NR2}. Bonacich \cite{Bo} introduces a family of measures tailored to a variety of
applications, De la Cruz Cabrera et al. \cite{DCMR} use subgraph centrality to the nodes 
of the line graph associated with original graph (which yields an importance measure for 
the edges of original graph), Ghosh et al. \cite{GBS} measure the effective resistance of 
a graph, Massei and Tudisco \cite{MT} determine the robustness of a network with 
Krylov subspace methods. Methods that are based on functions of the adjacency matrix for
the graph are discussed by Schweitzer \cite{Sc} and in \cite{NR4}. These functions may be
the matrix exponential, but are not limited to this matrix function.} 
The edges chosen depend on how their importance 
is measured and the latter often depends on the application being considered. It is outside 
the scope of the present paper to compare all available methods. The approach of the
present paper distinguishes itself by being based on perturbation analysis of the Perron
and Fiedler vectors and the condition numbers of these vectors.

This paper is organized as follows. Section \ref{sec2} analyzes the sensitivity of the 
Perron and Fiedler values to edge weight perturbations, as well as the conditioning of the
Perron and Fiedler vectors. In Section \ref{sec3} we investigate the effect of targeted 
edge removals on the Perron and Fiedler eigenpairs, and analyze some illustrative examples
of both synthetic and real-world networks. Section \ref{sec4} is concerned with multiplex 
networks. Conclusions are drawn in Section \ref{sec5}.

\section{The effect of structural perturbations on the Perron and Fiedler values}
\label{sec2}
Throughout this paper the superscript $(\cdot)^H$ stands for transposition and complex 
conjugation, $\|\cdot\|_2$ denotes the Euclidean vector norm or the spectral matrix norm, 
and $\|\cdot\|_F$ is the Frobenius matrix
norm. We first review a classical 
result due to Wilkinson that gave rise to the notion of eigenvalue condition number as the
reciprocal of the cosine of the angle between the corresponding left and right 
eigenvectors; see \cite[Section 2]{Wi}. We state this result for a general square matrix 
$M\in\C^{n\times n}$. 

\begin{proposition}(\cite{Wi})\label{cond}
Let $M\in\C^{n \times n}$, and let $x$ and $y$ be right and left eigenvectors of unit 
Euclidean norm associated with a simple eigenvalue $\lambda$. Define the matrix
$M^{(\varepsilon)}=M+\varepsilon E$, where $\varepsilon>0$ is small and 
$E\in{\mathbb C}^{n\times n}$ is a matrix with $\|E\|_2=1$. Let 
$\lambda^{(\varepsilon)}=\lambda+\delta\lambda$ be the eigenvalue of $M^{(\varepsilon)}$ 
that corresponds to $\lambda$, i.e., there is a continuous mapping $t\to \lambda^{(t)}$ for 
$0\leq t\leq\varepsilon$ such that $\lambda^{(t)}=\lambda$ for $t=0$ and 
$\lambda^{(t)}=\lambda^{(\varepsilon)}$ for $t=\varepsilon$. Then 
\begin{equation}\label{delta}
\delta\lambda=\varepsilon \frac{y^HE x}{y^Hx}+{\mathcal O}(\varepsilon^2) \quad\mbox{as}
\quad \varepsilon\searrow 0
\end{equation}
and
\begin{equation*}
\left |\frac{y^HEx}{y^Hx}\right| \leq\frac{\|y\|_2\|E\|_2\|x\|_2}{|y^Hx|}=
\frac{1}{\cos\theta_{x,y}}
\end{equation*}
with $\cos\theta_{x,y} := |y^Hx|$. The upper bound is attained for $E=yx^H$. 
The spectral norm of $E$ may be replaced by the Frobenius norm.
\end{proposition}

The condition number $\kappa(\lambda)$ of the simple eigenvalue $\lambda$ of a matrix 
$M\in\C^{n \times n}$ measures the worst-case effect on $\lambda$ of a tiny perturbation
of $M$. Proposition \ref{cond} suggests the following definition of the condition number; 
see \cite{Wi}. 

\begin{definition}
The condition number of a simple eigenvalue $\lambda$ of a matrix $M\in\C^{n\times n}$ is
defined by
\[
\kappa(\lambda)=\frac{1}{\cos\theta_{x,y}}
\]
with $\cos\theta_{x,y}$ given in Proposition \ref{cond}. The matrix
\begin{equation}\label{W}
W=yx^H.
\end{equation}
with the vectors $x$ and $y$ as in Proposition \ref{cond} is referred to as the 
\emph{Wilkinson perturbation} for $\lambda$.
\end{definition}

\begin{remark} \label{rmk1}
Clearly, $\kappa(\lambda)\geq 1$. When the matrix $M$ is normal, the right and left
eigenvectors can be chosen to coincide. Then $\kappa(\lambda)=1$ { and  $W=xx^H$}.
\end{remark}

We also are interested in the effect of {    structure-preserving perturbations, i.e., in
perturbation matrices whose non-trivial entries correspond to non-trivial entries of the 
matrix under consideration}; see, e.g., \cite{NP,KKT}. To investigate the structured 
eigenvalue sensitivity, we make use of the following result.

\begin{proposition}(\cite{NP})\label{strcon}
Let $M\in\mathbb{C}^{n \times n}$, and let $x$ and $y$ be right and left eigenvectors of 
unit Euclidean norm associated with the simple eigenvalue $\lambda$. Let 
${\cal S}\subseteq \mathbb{C}^{n\times n}$ denote the subspace of the complex matrices 
with the same sparsity structure as $M$, and let $N|_{\cal S}$ be the  projection of 
the matrix $N\in\mathbb{C}^{n \times n}$ onto ${\cal S}$.  Define the matrix
$M^{(\varepsilon)}=M+\varepsilon E$, where $\varepsilon>0$ is small and
$E\in {\cal S}$ is a matrix with $\|E\|_F=1$. Let 
$\lambda^{(\varepsilon)}=\lambda+\delta\lambda$ be the eigenvalue of $M^{(\varepsilon)}$ 
that corresponds to $\lambda$. Then \eqref{delta} holds and
\begin{equation*}
\left |\frac{y^HEx}{y^Hx}\right |\leq \frac{\| (yx^H)|_{\cal S} \|_F}{\cos\theta_{x,y}},
\end{equation*}
with $\cos\theta_{x,y} := |y^Hx|$. The upper bound is attained for 
$E=\frac { (yx^H)|_{\cal S} }{ \|(yx^H)|_{\cal S}\|_F}$. 
\end{proposition}

Proposition \ref{strcon} suggests the following definition of the ${\cal S}$-structured
condition number.

\begin{definition}(\cite{NP})
The ${\cal S}$-structured condition number of a simple eigenvalue $\lambda$ of a matrix 
$M\in\mathbb{C}^{n\times n}$ with unit right and left eigenvectors $x$ and $y$,
respectively, is defined by
$$
\kappa_{\cal S}(\lambda)=\kappa(\lambda)\|(yx^H)|_{\cal S}\|_F
$$
and the ${\cal S}$-structured analogue of the Wilkinson perturbation for $\lambda$ is
$$
W_{\cal S}=\frac { (yx^H)|_{\cal S} }{ \|(yx^H)|_{\cal S}\|_F}.
$$
\end{definition}

\begin{remark} 
Since $1=\|y\|_{{2}}\|x\|_{{2}} = \|yx^H\|_F \geq \|(yx^H)|_{\cal S}\|_F$, we have 
$\kappa_{\cal S}(\lambda)\leq \kappa(\lambda)$.
\end{remark}

\begin{remark} \label{rmk4}
{    When the matrix $M$ is normal, $\kappa_{\cal S}(\lambda)=\|(xx^H)|_{\cal S}\|_F
$ and $
W_{\cal S}=\frac { (xx^H)|_{\cal S} }{ \|(xx^H)|_{\cal S}\|_F}
$;  see Remark \ref{rmk1}.}
\end{remark}

\subsection{The impact of an edge removal on the Perron value}\label{subsec21}
Let the entries $a_{ij}=a_{ji}$, with $i\ne j$, of the symmetric adjacency matrix $A$ for
the graph ${\cal G}$ be positive. Then there is an edge between nodes $i$ and $j$, which 
we denote by $e(i\leftrightarrow j)$. To investigate the effect on the robustness of the 
network ${\cal G}$ when decreasing the weights $a_{ij}=a_{ji}$ of this edge, we 
consider the perturbed adjacency matrix 
\begin{equation}\label{Atau}
A^{(\tau)}=A-\tau a_{ij}(e_ie_j^T+e_je_i^T), \quad 0<\tau\leq1,
\end{equation} 
where $e_i=[0,\ldots,0,1,0,\ldots,0]^T$ denotes the $i$th vector of the canonical basis of
$\R^n$. Letting $\tau=1$ corresponds to removing the edge $e(i\leftrightarrow j)$.

Consider the Perron eigenpair $(\rho,u)$ for the matrix $A$. All entries of $A$ are 
nonnegative, and since the graph ${\cal G}$ is connected, the adjacency matrix $A$ is 
irreducible. Therefore, the Perron-Frobenius Theorem applies. It states that $\rho$ is a 
simple positive eigenvalue that decreases when any positive entry $a_{ij}$ of $A$ is 
decreased; see, e.g., \cite[Section 2]{PF}. 

Let $i\ne j$ and apply Proposition \ref{cond} with 
\begin{equation}\label{appl}
M=A,\qquad E=\frac{e_ie_j^T+e_je_i^T}{\sqrt{2}},\qquad \varepsilon=-\sqrt{2}\,\tau a_{ij},
\end{equation}
where we note that $\|E\|_F=1$. The perturbed adjacency matrix $A^{(\tau)}$ in 
\eqref{Atau} is symmetric and nonnegative with all diagonal entries zero {    due to the 
absence of self-loops in the graph}.  We are interested in determining the impact on $\rho$ 
of {    changes of} the weights of the edge $e(i\leftrightarrow j)$ of ${\cal G}$. Let 
$\rho^{(\tau)}$ denote the spectral radius of $A^{(\tau)}$. {\  Assuming that $i>j$, 
we obtain from \eqref{delta} and \eqref{appl} that the impact of edge $e(i\leftrightarrow j)$ 
on $\rho$ can be estimated as follows.
\begin{equation} \label{srho}
\frac{\rho-\rho^{(\tau)}}{\rho}\approx\frac{\tau a_{ij} }{\rho} \frac{u^T(e_ie_j^T+e_je_i^T)u}{u^Tu}=
\frac{2\tau}{\rho} a_{ij} u_iu_j,\qquad i>j.
\end{equation}}

Since the matrix $A$ is symmetric, the right and left Perron vectors coincide. The  
eigenvalue condition number of $\rho$ therefore is $1$ and \cite[eqs. (2)-(3)]{NR1} 
simplify. The Wilkinson perturbation \eqref{W} associated with $\rho$ is the nonnegative
symmetric matrix $W=[u_iu_j]\in\R^{n\times n}$, whose $(i,j)$th entry appears in the 
right-hand side of \eqref{srho}. 

Let $R^{(\rho)}\in\R^{n\times n}$ be the strictly lower triangular matrix, whose $(i,j)$th
entry for $i>j$ approximates the impact on the Perron { value} of the removal of the edge 
$e(i\leftrightarrow j)$ from the graph ${\cal G}$ according to the right-hand side of 
\eqref{srho}. We will refer to $R^{(\rho)}$ as the {\it Perron impact matrix}. 

Let ${\cal L}\subseteq\R^{n\times n}$ denote the cone of the real nonnegative $n\times n$
matrices with the same sparsity structure as the strictly lower triangular portion of 
$A$ and let $N|_{\cal L}$ denote the projection of a nonnegative matrix 
$N\in\R^{n\times n}$ onto ${\cal L}$. 

\begin{proposition}
The Perron impact matrix is given by 
\begin{equation}\label{Rmat}
R^{(\rho)}=\frac{2}{\rho} A|_{\cal L}\circ W|_{\cal L},
\end{equation}
where $\circ$ denotes the Hadamard product. Moreover,
\begin{equation}\label{Rmatbd}
\|R^{(\rho)}\|_F\leq\frac{\sqrt 2}{\rho} \kappa_{\cal S}(\rho)\|A|_{\cal L}\|_F,
\end{equation}
where ${\cal S}\subseteq\mathbb{R}^{n\times n}$ is the cone of all nonnegative { symmetric} matrices 
with the same sparsity structure as $A$. 
\end{proposition}

\proof
It follows from \eqref{srho} that the Perron impact matrix can be determined by multiplying
entry-wise the strictly lower 
triangular portion of $W$ and the strictly lower triangular portion of $A$ multiplied by 
the factor $\frac{2}{\rho}$. The right-hand side is nonvanishing if and only if 
$a_{ij}>0$, which corresponds to an edge $e(i\rightarrow j)$. This shows \eqref{Rmat}. To
verify \eqref{Rmatbd}, we note that the Frobenius norm is Hadamard submultiplicative, 
i.e.,
\[
\|M\circ N\|_F\leq \|M\|_F \|N\|_F,\qquad\forall\ M,N\in\C^{n\times n};
\]
see, e.g., \cite[p. 335]{HJ2}. Therefore,
\[
\|R^{(\rho)}\|_F\leq\frac{2}{\rho}\|A|_{\cal L}\|_F\|W|_{\cal L}\|_F.
\]
The proof now follows by observing that  
\[
\|W|_{\cal L}\|_F=\frac{\sqrt 2}{2}\kappa_{\cal S}(\rho),
\]
because $\|W|_{\cal S}\|_F^2=2\|W|_{\cal L}\|_F^2$ and 
$\kappa_{\cal S}(\rho)= \|W|_{\cal S}\|_F$.
\endproof

\begin{remark}
If the network is unweighted, then $A|_{\cal L}\circ W|_{\cal L}=W|_{\cal L}$. This gives 
\[
\|R^{(\rho)}\|_F=\frac{\sqrt {2}}{\rho} \kappa_{\cal S}(\rho).
\]
\end{remark}

\subsection{The impact of an edge removal on the Fiedler value}\label{subsec22}
Assume that the Fiedler eigenpair $(\mu,v)$ is simple. Let the matrix $A^{(\tau)}$ be 
defined by \eqref{Atau} and let $D^{(\tau)}=\mathrm{diag}[A^{(\tau)} 1_n]$. The graph 
Laplacian determined by the graph associated with $A^{(\tau)}$ is given by
\begin{eqnarray*}
L^{(\tau)}&=&D^{(\tau)}-A^{(\tau)}=(D-\tau a_{ij}(e_ie_i^T+e_je_j^T))-
(A- \tau a_{ij}(e_ie_j^T+e_je_i^T)) \\
&=& L - \tau a_{ij}((e_ie_i^T+e_je_j^T)-(e_ie_j^T+e_je_i^T)), \qquad 0<\tau\leq 1.
\end{eqnarray*}
Let $\mu^{(\tau)}$ denote the second smallest eigenvalue of $L^{(\tau)}$. We apply 
Proposition \ref{cond} with
\[
M=L=D-A,\quad M^{(\varepsilon)}=L^{(\tau)}=D^{(\tau)}-A^{(\tau)},
\]
\begin{equation}\label{Emat}
E=\frac{e_ie_i^T+e_je_j^T-e_ie_j^T-e_je_i^T}{2},
\end{equation}
and $\varepsilon=-2\,\tau a_{ij}$. Note that $E$ is a symmetric positive semidefinite 
matrix of rank one with $\|E\|_F=1$. It follows from  Weil's eigenvalue perturbation 
inequalities (see, e.g., \cite[pp. 182--183]{HJ}) applied to a symmetric positive 
semidefinite perturbation of rank one of a symmetric matrix, that $\mu^{(\tau)}$ is smaller
than or equal to $\mu$, and that the third smallest  eigenvalue of $L^{(\tau)}$ is larger 
than or equal to $\mu$ and smaller than or equal to the third smallest eigenvalue of $L$.

\begin{remark} \label{decmu}
{    We recall that the Braess paradox \cite{Br} is the phenomenon showing that removing 
an edge may increase a measure of connectivity of a network rather than reduce it; see 
\cite{ABCMP} for a discussion. While removing an edge always decreases the algebraic 
connectivity $\mu$ of a network (cf. Proposition \ref{pro8} below), it is interesting to 
note that, conversely, the following statement may be false:
$$ \mathrm{If }\;  0<\tau_1<\tau_2<1, \; \mathrm{ then } \; \mu^{(\tau_1)} > \mu^{(\tau_2)}. $$ 
}
\end{remark}

{\  We are interested in estimating the impact on the Fiedler value $\mu$ of reducing 
the weights of the edge $e(i\leftrightarrow j)$ that determine the matrix \eqref{Emat}. 
From \eqref{delta} we obtain the following first-order estimate for the impact of edge 
$e(i \leftrightarrow j)$ on $\mu$:
\begin{equation}\label{smu}
\frac{\mu-\mu^{(\tau)}}{\mu}\approx\frac{\tau a_{ij}}{\mu}\frac{v^T (e_ie_i^T+e_je_j^T-e_ie_j^T-
e_je_i^T)v}{v^Tv}=\frac{\tau a_{ij}}{\mu} (v_i^2+ v_j^2-2v_iv_j)=\frac{\tau}{\mu} 
a_{ij} (v_i-v_j)^2.
\end{equation}
In particular, we are
interested in removing this edge. This corresponds to $\tau=1$. }

\begin{proposition}
Let $Y=[y_{ij}]\in\R^{n\times n}$ be the strictly lower triangular matrix with entries 
$y_{ij}=(v_i-v_j)^2$, $i>j$, and let $Y|_{\cal L}$ denote the projection of $Y$ onto the 
sparsity structure of the matrix $ A$. Define the \emph{Fiedler impact matrix}
\begin{equation}\label{Rmu}
R^{(\mu)}=\frac{1}{\mu} A|_{\cal L}\circ Y|_{\cal L}.
\end{equation}
Then the $(i,j)$th entry of $R^{(\mu)}$ for $i>j$ approximates the impact of the edge
$e(i\leftrightarrow j)$ in ${\cal G}$ on the Fiedler value. When the graph ${\cal G}$ is
unweighted, 
\begin{equation}\label{Rmubd}
\|R^{(\mu)}\|_F=\frac{1}{\mu} \|Y|_{\cal L}\|_F\leq \frac{{ 2} \, \sqrt{m}}{\mu}.
\end{equation}
\end{proposition}

\proof
The relation \eqref{Rmu} follows from the right-hand side of \eqref{smu}. We turn to 
\eqref{Rmubd}. When the graph is unweighted, all nonvanishing entries of $A$ are one.
Therefore, $A|_{\cal L}\circ Y|_{\cal L}=Y|_{\cal L}$. Since 
$$
 \max_{v_i^2+v_j^2 \le 1} (v_i-v_j)^2 =2
$$
is attained for $v_i=\pm\frac{\sqrt{2}}{2}$ and $v_j=\mp\frac{\sqrt{2}}{2}$,
we obtain $\|Y|_{\cal L}\|_F\leq { 2} \, \sqrt{m}$.
\endproof

\subsection{Sensitivity of the Perron and Fiedler vectors}\label{subsec23}
We are interested in how an eigenvector $x$ of a matrix $M\in\C^{n\times n}$ associated 
with a simple eigenvalue is affected by a perturbation of $M$. The condition number of an
eigenvector $x$ provides a bound for the angle between $x$ and the corresponding perturbed
eigenvector $x^{(\varepsilon)}$. {    We will apply this result below to derive a lower
bound for the condition number of the Perron vector of the adjacency matrix $A$ for the
graph ${\cal G}$.}

\begin{definition}(\cite{St01})\label{condvec1}
Let $M\in\mathbb{C}^{n \times n}$ and let $x$ be an eigenvector of unit Euclidean norm 
associated with a simple eigenvalue $\lambda$. Assume that the columns of the matrix 
$U\in\C^{n\times (n-1)}$ form an orthonormal basis for $\text{Range}(M-\lambda I_n)$, 
where $I_k$ denotes the identity matrix in $\C^{k\times k}$. The condition number of $x$ 
(i.e., the condition number of the one-dimensional invariant subspace spanned by $x$) is 
defined as 
\begin{equation}\label{kappadef}
\kappa(x)=\|(\lambda I_{n-1} - U^HMU)^{-1}\|_2^{-1}.
\end{equation}
\end{definition}

Let $x$ be an eigenvector of unit Euclidean norm associated with the simple eigenvalue 
$\lambda$ of the matrix $M\in\C^{n \times n}$. Consider the matrix 
$M^{(\varepsilon)}=M+\varepsilon E$, where $\varepsilon>0$ is small and 
$E\in\C^{n\times n}$ is a matrix with $\|E\|_F=1$. Let $x^{(\varepsilon)}$ be the unit 
eigenvector of $M^{(\varepsilon)}$ that corresponds to $x$, i.e., there is a continuous 
mapping $t\to x^{(t)}$ for $0\leq t\leq\varepsilon$ such that $x^{(t)}=x$ for $t=0$ and 
$x^{(t)}=x^{(\varepsilon)}$ for $t=\varepsilon$. Then 
\[
\sin \theta_{x,x^{(\varepsilon)}}\leq \kappa(x)\varepsilon,
\]
where
\[
\sin \theta_{x,x^{(\varepsilon)}} := \sqrt{1-\cos^2 \theta_{x,x^{(\varepsilon)}}},\quad  
\cos \theta_{x,x^{(\varepsilon)}} :=|x^H x^{(\varepsilon)}|,
\]
and $\kappa(x)$ is given by \eqref{kappadef}; see Stewart \cite[pp. 48--50]{St01} for 
details.

Consider the special case when the matrix $M\in\C^{n\times n}$ is normal and denote its 
spectrum by $\Lambda(M)$. Following the argument in \cite{NR}, regard the the reciprocal 
of the expression in the right-hand side of \eqref{kappadef}. It is straightforward to 
show that the upper bound
\[
\|(\lambda I_{n-1} - U^HMU)^{-1}\|_2 \leq \min_{{\substack{\gamma\neq\lambda\\
\gamma\in\Lambda(M)}}} |\lambda-\gamma|
\]
is attained because $M$ is unitarily diagonalizable. Therefore, for a normal matrix $M$
the condition number of a unit eigenvector $x$ only depends on how well the associated 
eigenvalue $\mu$ is separated from the other eigenvalues of the matrix. This observation 
leads to the following proposition, which is shown in \cite[p. 51]{St01}.

\begin{proposition}\label{codvecdef}
Let $M\in\C^{n \times n}$ be a normal matrix and let $x$ be a unit eigenvector associated 
with a simple eigenvalue $\lambda$. The condition number of $x$ (i.e., the condition 
number of the one-dimensional invariant subspace spanned by $x$) is given by
\[
\kappa(x)=\frac{1}{\displaystyle{\min_{\substack{\gamma\neq\lambda\\ \gamma\in\Lambda(M)}} 
|\lambda-\gamma|}}.
\]
\end{proposition}

Let $A\in\R^{n\times n}$ be the symmetric irreducible adjacency matrix associated with an 
undirected and connected graph $\cal G$. The eigenvalue $\gamma^*$ of $A$ that is closest
to the Perron value $\rho$, i.e., 
\[
\gamma^*=\arg\min_{\substack{\gamma\neq\rho\\ \gamma\in\Lambda(A)}}|\rho-\gamma|,
\]
satisfies the following bounds that can be used in the analysis of the condition number 
of the Perron vector.

\begin{proposition} \label{pgam}
Let $\rho$ denote the Perron value of the symmetric irreducible adjacency matrix
$A\in\R^{n\times n}$. Then
\[
-\frac{\rho}{n-1}\leq\gamma^*<\rho.
\]
\end{proposition}

\proof
Let the matrix $A$ have the eigenvalues 
$\lambda_1\leq\lambda_2\leq \ldots \leq \lambda_{n-1}\leq \lambda_n$. Due to the 
Perron-Frobenius Theorem, $\lambda_n=\rho>\max\{\lambda_{n-1},0\}$. Moreover, 
since $A$ has zero trace, one has
\[
0=\sum_{i=1}^n \lambda_i=\sum_{i=1}^{n-1}\lambda_i+\rho\leq (n-1)\lambda_{n-1}+\rho.
\]
The proof now follows from $\gamma^*=\lambda_{n-1}$.
\endproof

\begin{proposition}\label{prop7}
Let $A\in\R^{n\times n}$ be a symmetric irreducible adjacency matrix. Let $\rho$ denote 
the Perron value of $A$ and let $u$ be the Perron vector. Then
\begin{equation}\label{lwrbd}
\kappa(u)\geq \frac{n-1}{n \rho}.
\end{equation}
\end{proposition}

\proof
It follows from Proposition \ref{pgam} that
\[
\rho - \gamma^*\leq \frac{\rho}{n-1}+\rho= \rho\;\frac{n}{n-1}.
\]
Therefore, we obtain from Proposition \ref{codvecdef}, 
\[
\kappa(u)=\frac{1}{{ \rho}-\gamma^*} \geq \frac{n-1}{n \rho},
\]
which concludes the proof. \endproof

\begin{remark}
The Perron-Frobenius Theorem \cite[Section 2]{PF} assures that decreasing the weight of 
any edge of ${\cal G}$ decreases the Perron value $\rho$. The preceding result shows that 
if $\rho$ decreases, then the lower bound \eqref{lwrbd} for the condition number of the 
Perron vector $u$ increases.
\end{remark}

The analysis of the condition number of the Fiedler vector $v$, assuming that the Fiedler 
value $\mu$ is simple, is much more difficult. A reason for this is that the Fiedler value
might not change monotonically with a decrease of edge weights {    (cf. Remark \ref{decmu})}.

\section{Perturbations of the Perron and Fiedler eigenpairs induced by edge removal}
\label{sec3}
In many applications of graphs, the edge weights are estimated. They therefore are subject
to error. We are interested in analyzing how sensitive the Perron and Fiedler vectors are 
to perturbations of the edge weights. 

Consider the removal of the edge $e(i\leftrightarrow j)$ from a graph ${\cal G}$ and 
assume that $i>j$. Then the $(i,j)$th entries of the matrices $R^{(\rho)}$ and $R^{(\mu)}$ are 
important. The former entry is given by $\frac{2}{\rho}a_{ij} u_iu_j$; we refer to it as 
the \emph{relative Perron edge importance} for $e(i \leftrightarrow j)$. The $(i,j)$th 
entry of $R^{(\mu)}$ is given by $\frac{1}{\mu}a_{ij} (v_i-v_j)^2$ and is referred to as 
the \emph{relative Fiedler edge importance}. The larger the relative Perron edge 
importance, the more will a small reduction of the edge weights $a_{ij}=a_{ji}$ reduce the 
Perron value $\rho$. This suggests that in order to reduce the Perron value, it may be 
beneficial to remove an edge $e(i\leftrightarrow j)$ with a large relative Perron edge 
importance. Similarly, if an edge $e(i\leftrightarrow j)$ { that} is associated with a relatively 
large Fiedler edge importance is removed, then there is a risk that the reduced network 
is not connected.


Since the Perron value $\rho$ is a 
simple eigenvalue, it follows from Proposition \ref{codvecdef} with $x=u$ and 
$\lambda=\rho$ that the Perron vector $u$ is sensitive to perturbations of the adjacency
matrix $A$ when $\rho$ and the closest eigenvalue, $\lambda_{n-1}$, are close, because 
then the condition number $\kappa(u)$ is large; cf. the proof of Proposition \ref{prop7}. 
A large value of $\kappa(u)$ indicates that the Perron vector may be quite sensitive to 
perturbations of the weights of the graph. As mentioned above, the weights may be subject 
to errors. Therefore, the edge removal suggested by the entries of the Perron vector might
not be appropriate when $\kappa(u)$ is large. Moreover, removing an edge from ${\cal G}$
induces a perturbation in the adjacency matrix $A$. 
When $\kappa(u)$ is large, it may be
important to compute the Perron vector for the reduced adjacency matrix obtained by
the edge removal before seeking to remove the next edge. 
{ When proceeding with the bipartition of the network, as illustrated in the following 
subsection, it may be beneficial to check the sensitivity of the Perron vector for each 
adjacency matrix obtained by the edge removal before seeking to remove the next edge.}

We turn to the Fiedler sensitivity analysis. Let $e(i\leftrightarrow j)$ be an edge of
${\cal G}$. We are interested in how the spectrum of the graph Laplacian $L$ associated 
with $\cal G$ changes when removing this edge from the graph. 
The following proposition is shown in \cite[Sect. 3.2]{BH}.

\begin{proposition}(\cite{BH})\label{pro8}
Let $0=\alpha_1\leq\alpha_2\leq\ldots\leq\alpha_{n-1}\leq \alpha_n$ be the 
eigenvalues of the graph Laplacian $L$ associated with {    an undirected weighted graph} 
${\cal G}$ and let $0=\widetilde{\alpha}_1\leq \widetilde{\alpha}_2\leq\ldots\leq\widetilde{\alpha}_{n-1}\leq
\widetilde{\alpha}_n$ denote the eigenvalues of the graph Laplacian $\widetilde{L}$ 
associated with the graph $\widetilde{\cal G}$ obtained by removing the edge 
$e(i\leftrightarrow j)$ from $\cal G$. Then the eigenvalues $\widetilde{\alpha}_j$ and 
$\alpha_k$ interlace, i.e.,
\begin{equation}\label{incaps}
0=\widetilde{\alpha}_1=\alpha_{1} \leq\widetilde{\alpha}_2\leq\alpha_{2}\leq\dots\leq
\widetilde{\alpha}_{n-1} \leq \alpha_{n-1}\leq  \widetilde{\alpha}_n \leq \alpha_{n}.
\end{equation}
\end{proposition}

If the Fiedler eigenpair is simple, then the same considerations made above for the
Perron eigenpair apply. In this case $\kappa(v)$ is the reciprocal of the distance 
between $\mu$ and the closest eigenvalue of the graph Laplacian.

Numerical tests with both synthetic and real-world complex networks are reported in the
following subsection. Our analysis of synthetic graphs accounts for bipartivity \cite[p. 11]{Bi} 
(cf. Examples \ref{ex1} and \ref{ex2}) and the presence of anti-communities, that is the
presence of subsets of nodes that are loosely connected among each other but highly 
connected to the rest of the network; see, e.g., \cite{CNRR,FT} (cf. Example \ref{ex3}). 

\begin{remark}
We recall that a network is said to be \emph{bipartite} if the set of nodes can be 
subdivided into two nonempty subsets such that there are no edges between nodes in the 
same subset. Bipartivity is an important structural property. 

A sequence of $k$ edges and $k+1$ nodes such that 
\[
\{e(i_1\leftrightarrow i_2), e(i_2\leftrightarrow i_3),\ldots,e(i_k\leftrightarrow 
i_{k+1})\}
\]
are said to form a \emph{walk}. The length of a walk is the sum of the weights of the 
edges in the walk, i.e., $\sum_{j=1}^k a_{i_j,i_{j+1}}$. If the edges in a walk are 
distinct, then the walk is referred to as a \emph{path}. A graph is bipartite if there are
no closed walks with an odd number of edges.

It follows from the application of Gershgorin disks that all eigenvalues of the 
\emph{normalized graph Laplacian}
\begin{equation}\label{normL}
D^{-\frac{1}{2}}LD^{-\frac{1}{2}}=I_n-D^{-\frac{1}{2}}AD^{-\frac{1}{2}}
\end{equation}
live in the interval $[0,2]$. The multiplicity of the eigenvalue zero of \eqref{normL} 
equals the number of 
connected components of ${\cal G}$; the multiplicity of the eigenvalue $2$ is the number 
of bipartite components of ${\cal G}$ with two or more nodes; see \cite{BC,Ch}. The latter
multiplicity coincides with the multiplicity of the eigenvalue $-1$ in the spectrum of the
matrix $A_D:=D^{-\frac{1}{2}}AD^{-\frac{1}{2}}$. This matrix is known as the 
\emph{reduced adjacency matrix}. Determining the multiplicity of the eigenvalue $-1$ 
of the matrix $A_D$ is an effective way to determine whether a graph is bipartite. 
\end{remark}

\subsection{Numerical examples}
For small to moderately sized networks, Perron and Fiedler eigenpairs can easily be 
evaluated by using the MATLAB function \texttt{eig}. For large-scale symmetric adjacency 
matrices, these quantities can be computed by the MATLAB function \texttt{eigs}. This
function uses a shift-and-invert strategy for the evaluation of the desired eigenpair and 
therefore requires the solution of linear systems of equations with shifted adjacency 
matrices $A-\alpha_j I_n$, $\alpha_j\in\R$. Alternatively, the IRBL algorithm 
\cite{BCR1,BCR2}, which implements a restarted Lanczos method for computing extreme 
eigenpairs of a large symmetric matrix without using a shift-and-invert strategy, can be 
applied. This algorithm does not require the solution of large linear systems of equations 
and, therefore, can be applied to very large adjacency matrices. The numerical tests 
reported in this paper have been carried out using MATLAB R2024b on a 3.2 GHz {    Intel 
Core i7 iMac (6 cores). All computations are performed with about $15$ significant
decimal digits}. 

The first three examples discuss synthetic connected undirected networks with $n=6$
nodes. Let the eigenvalues of the associated adjacency matrices be sorted in non-decreasing 
order, 
\[
\lambda_1\leq \ldots \leq \lambda_{5}<\lambda_{6}.
\]
Also sort the eigenvalues of the corresponding graph Laplacians in non-decreasing order,
\[
0=\alpha_1\leq\alpha_2\leq\alpha_3\leq \ldots \leq\alpha_{6}.
\]

\begin{example}\label{ex1}
\begin{figure}
\centering
\begin{tabular}{cc}
\includegraphics[scale = 0.4]{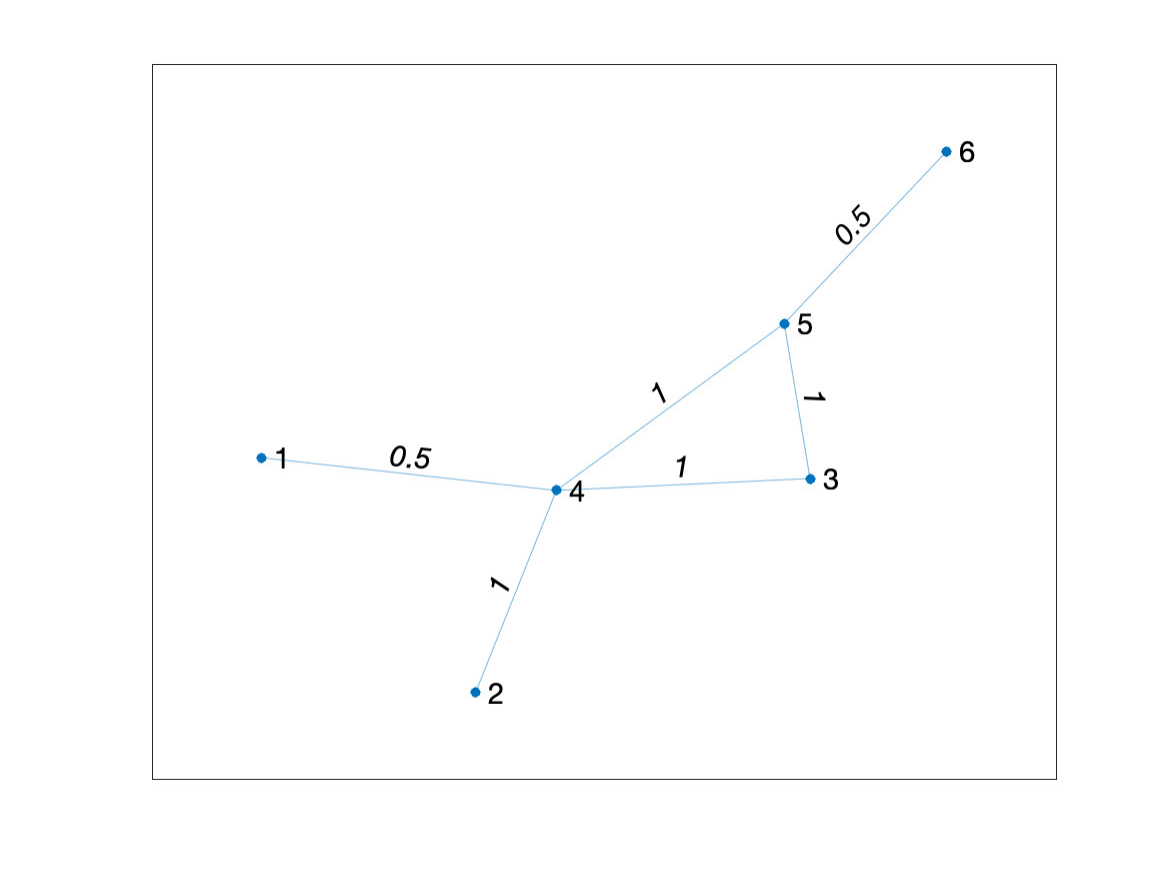} 
\includegraphics[scale = 0.4]{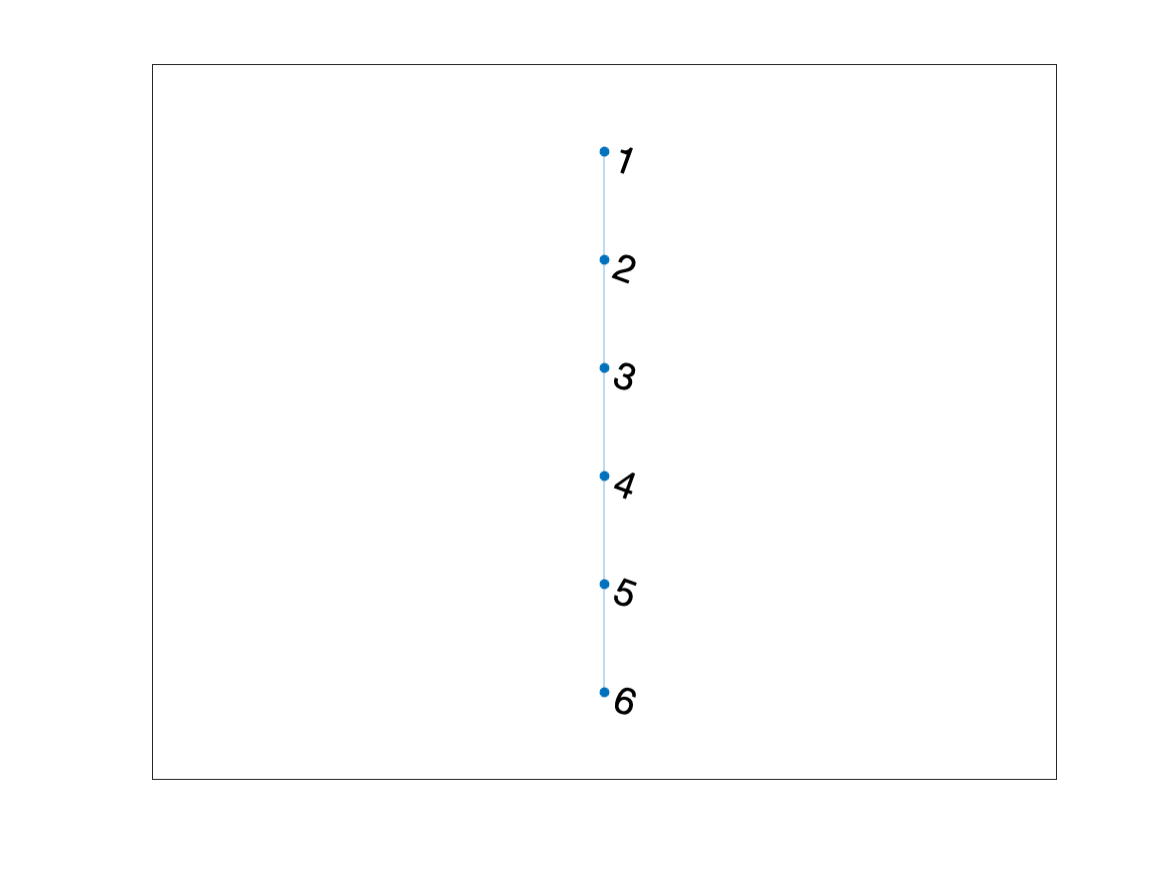} 
\end{tabular}
\caption{Graphs considered in Example \ref{ex1} (left) and Example \ref{ex2} (right).}
\label{FIG1}
\end{figure}

Consider the graph $\cal G$ depicted in Figure \ref{FIG1} (left); the edge weights are
marked in the figure. The spectral radius of 
the associated adjacency matrix $A$, which is irreducible, is $\rho=\lambda_{6}=2.2431$. 
The second smallest eigenvalue of the graph Laplacian $L$ associated with $\cal G$ is 
$\mu=\alpha_2=0.4100$. This eigenvalue is simple, i.e., $0<\alpha_2<\alpha_3$.
The Perron and Fiedler vectors are 
\begin{eqnarray}
\nonumber
u&=& \left[0.1355, \, 0.2710, \, 0.5034, \, 0.6079, \, 0.5213,\, 0.1162\right]^T,\\
\label{vvv}
v&=& \left[-0.6083, \, -0.1855, \, 0.0162, \, -0.1094, \, 0.1352,\, 0.7517\right]^T.
\end{eqnarray}
The second largest eigenvalue of $A$ is $\lambda_5=0.4574$ and the third smallest 
eigenvalue of $L$ is $\alpha_3=0.6436$. Therefore,
\[
\kappa(u)=\frac{1}{\rho-\lambda_5}=0.5600, \qquad 
\kappa(v)=\frac{1}{\min\{\alpha_3-\mu,\mu\}}=4.2808.
\]

The edges whose removal have the largest impact on $\rho$ are $e(5 \leftrightarrow 4)$ 
and $e(4 \leftrightarrow 3)$; both edges have weight one. The relevant entries of the 
Perron impact matrix $R^{(\rho)}$ are $0.2826$ and $0.2728$, respectively. 

The edges that have the largest impact on $\mu$ are $e(6 \leftrightarrow 5)$ and 
$e(4 \leftrightarrow 1)$; both edges have weight $0.5$. The relevant entries of the 
Fiedler impact matrix $R^{(\mu)}$ are $0.4634$ and $0.3034$, respectively. 

Table \ref{T_1} shows several eigenvalues associated with the adjacency matrices and graph 
Laplacians for the graph ${\cal G}$ and for some graphs that are obtained from ${\cal G}$ 
by removing selected edges. In detail, the table reports the Perron value and the two 
closest eigenvalues of the adjacency matrix $A$, as well as the condition number of the 
Perron vector $u$ {    as computed with MATLAB}. Table \ref{T_1} also shows eigenvalues of
the graph Laplacian $L$ and the condition number of the Fiedler vector $v$. 
 
The graphs ${\cal G}_j$, $j=1,2,\ldots,5$, are obtained from $\cal G$ as follows. We 
obtain graph ${\cal G}_1$ by removing the edge $e(5\leftrightarrow 4)$ from $\cal G$. This
is the edge whose removal has the largest impact on the spectral radius $\rho$. The 
resulting graph is connected. Graph ${\cal G}_2$ is obtained from $\cal G$ by removing the 
edge $e(4\leftrightarrow 3)$. This edge has the second largest impact on $\rho$. The graph
obtained is connected and has a slightly smaller condition number of the Perron vector 
than ${\cal G}_1$. We can see that the Perron value for ${\cal G}_1$ is smaller than
for ${\cal G}_2$ as can be expect from our analysis in Subsection \ref{subsec21}.

We obtain graph ${\cal G}_3$ by removing the edge $e(6\leftrightarrow 5)$ from $\cal G$. 
This edge has the largest impact on the algebraic connectivity $\mu$; the resulting graph 
is not connected. Therefore, the associated adjacency matrix is not irreducible (in fact, 
both vectors $u$ and $v$ have a zero entry). Moreover, the second smallest eigenvalue of 
the graph Laplacian for ${\cal G}_3$ is zero. We remark that the multiplicity of the zero 
eigenvalue is equal to the number of disconnected components of the graph. Graph 
${\cal G}_4$ is obtained from $\cal G$ by removing the edge $e(4\leftrightarrow 1)$. This 
is the edge with the second largest impact on $\mu$; the graph ${\cal G}_4$ is not 
connected.

Observe that the spectra of the graph Laplacians for ${\cal G}_1$, ${\cal G}_2$, 
${\cal G}_3$, and ${\cal G}_4$ satisfy the inequalities \eqref{incaps} when compared to
the spectrum of the graph Laplacian for ${\cal G}$. 

We now remove two edges from $\cal G$. Thus, we get graph ${\cal G}_5$ from $\cal G$ by 
removing the edges $e(5\leftrightarrow 4)$ and $e(4\leftrightarrow 3)$. These are the 
edges with the largest impact on $\rho$. The resulting graph is not connected, and neither
its spectral radius nor its algebraic connectivity are simple eigenvalues.

The spectrum of the adjacency matrix for ${\cal G}_5$ is symmetric with respect to the
origin. The symmetry indicates that the graph is made up of two bipartite components, 
which we refer to as ${\cal G}_5^{(1)}$ and ${\cal G}_5^{(2)}$. These components are 
disconnected. Interestingly, one can observe that the graphs ${\cal G}_5^{(1)}$ and 
${\cal G}_5^{(2)}$ correspond to the Fiedler bipartition suggested by the signs of 
\eqref{vvv}, with $N_1 = \{1,2,4\}$ and $N_2=\{3,5,6\}$. Both the graphs ${\cal G}_5^{(1)}$ and 
${\cal G}_5^{(2)}$ have the Perron value $\hat\rho=1.1180$ and the Fiedler value 
$\hat\mu=0.6340$. Since the eigenvalues preceding $\hat\rho$ and following $\hat\mu$ are
$0$ and $2.3660$, respectively, the condition numbers of the Perron and Fiedler vectors 
are 
\[
\kappa(u)=\frac{1}{\hat\rho}=0.8945, \qquad \kappa(v)=\frac{1}{\hat\mu}=1.5773.
\]

\begin{table}[h!tb]
\centering
\begin{tabular}{c | c c c c | c c c c }
&$\lambda_6$          & $\lambda_5$ & $\lambda_4$ &$\kappa(u)$  & $\alpha_2$ & $\alpha_3$& $\alpha_4$ & $\kappa(v)$\\
\hline
$\cal{G}$  & 2.2431 & $0.4574$ & $0.0000$&   $0.5600$  & { 0.4100} & 0.6436& 1.1475 & $4.2808$\ \\
${\cal{G}}_1$  &{ 1.6935} & $0.7950$ & $0.0000$&  $1.1130$ & { 0.2679} & 0.6340& 1.0000 & $3.7327$\\  
${\cal{G}}_2$  &1.7247 & $0.7247$& $0.0000$ &  $ 1.0000$ & { 0.3461} & 0.6340& 0.7378  & $3.4734$ \\
${\cal{G}}_3$  &2.2132 & $0.3592$ &$0.0000$ & $0.5394$  & 0.0000 & { 0.5635}& 1.0000&  $2.2910$ \\
${\cal{G}}_4$  &2.2020 & $0.4088$ &$0.0000$ &  $0.5577$  & 0.0000 & { 0.4902}& 1.1389 & $2.0400$ \\
${\cal{G}}_5$  &1.1180 &1.1180&  $0.0000$ &$-$ & 0.0000 & { 0.6340}& 0.6340  & $-$ \\

\end{tabular}
\caption{Example \ref{ex1}. Spectral measures relevant to graph ${\cal G}$ and derived graphs.} 
\label{T_1}
\end{table}
\end{example} 

\begin{example}\label{ex2}
Consider the unweighted graph ${\cal G}$ shown in Figure \ref{FIG1} (right) with $n=6$ 
nodes. Its adjacency matrix $A$ is the tridiagonal Toeplitz matrix with ones on the super- 
and sub-diagonals and zeros elsewhere. Both the eigenvalues and eigenvectors of $A$ are 
explicitly known; see, e.g., \cite{NR}. The graph ${\cal G}$ is bipartite and has spectral
radius $\rho=\lambda_{6}=1.8019$. The second smallest eigenvalue of the graph Laplacian 
matrix is the simple eigenvalue $\mu=\alpha_2=0.2679$. The Perron and Fiedler vectors are
\begin{eqnarray}
\nonumber
u&=& \left[0.2319,\, 0.4179,\, 0.5211,\, 0.5211,\, 0.4179,\, 0.2319\right]^T\\
\label{v}
v&=& \left[-0.5577, \, -0.4082, \, -0.1494, \, 0.1494, \, 0.4082,\, 0.5577\right]^T.
\end{eqnarray}
It follows that 
\[
\kappa(u)=\frac{1}{\rho-\lambda_5}=1.8019, \qquad \kappa(v)=\frac{1}{\mu}=3.7327.
\]
The edge whose removal has the largest impact on both $\rho$ and $\mu$ is 
$e(4\leftrightarrow 3)$. In fact, graph ${\cal G}_1$, which is obtained by removing this 
edge from $\cal G$, is not connected. It is made up of two bipartite components. In detail, 
all the eigenvalues of both the adjacency matrix and the graph Laplacian for ${\cal G}_1$ 
have algebraic multiplicity two, and the associated eigenvectors have three zero entries. 

The connected components of ${\cal G}_1$, which we denote by ${\cal G}_1^{(1)}$ and 
${\cal G}_1^{(2)}$, correspond to the Fiedler bipartition of ${\cal G}$ suggested by the
signs of the Fiedler vector in \eqref{v} with $N_1 = \{1,2,3\}$ and $N_2=\{4,5,6\}$. Both
components have the Perron value $\hat\rho=1.4142$ and the Fiedler value $\hat\mu= 1$.
Finally, the condition numbers of their Perron and Fiedler vectors are 
\[
\kappa(u)=\frac{1}{\hat\rho}=0.7071, \qquad \kappa(v)=\frac{1}{\hat\mu}=1.
\]
\end{example} 
\begin {remark}
Clearly, Example \ref{ex2} is not typical {    in the sense that, generically, the edge with
the largest impact on the spectral 
radius is not the same as the edge with the largest impact on the algebraic connectivity.}
For the sake of completeness, we report results of some numerical tests {    that
illustrate this}. In these tests we use random irreducible $n \times n$ adjacency matrices
for $n=6$, and take into account both cases with maximal and minimal numbers of edges 
(i.e.,  
$\frac{n(n-1)}{2}=15$ edges and $n-1=5$ edges, respectively). In $5\cdot 10^6$ tests with
adjacency matrices associated with complete graphs (i.e., graphs in which each node is 
connected to all the other nodes with a single edge), with weights uniformly distributed 
in the open interval $(0,1)$, it happened in about $2\cdot 10^4$ tests (less than $0.5\%$)
that the edge with the largest impact on the spectral radius was the same as the edge with 
the largest impact on the algebraic connectivity, whereas in $5 \cdot 10^6$ tests carried 
out with $6\times 6$ tridiagonal adjacency matrices with weights uniformly distributed in $(0,1)$, this
happened in about $5\cdot 10^4$ tests (about $1 \%$). 
\end{remark}

\begin{example} \label{ex3}
Consider the unweighted graph ${\cal G}$ shown in Figure \ref{FIG2} (left) with $n=6$ 
nodes. Its adjacency matrix has a zero leading square block of order $3$. The graph 
${\cal G}$ has an anti-community formed by the nodes $\{1, 2, 3\}$. The Perron value of 
${\cal G}$ is $\rho=\lambda_{6}=2.4142$ and the Perron vector is given by
\begin{equation*} 
u= \left[0.2209,\, 0.2209,\, 0.2209,\, 0.5334,\, 0.5334,\, 0.5334\right]^T.
\end{equation*}
{ The condition number of the Perron vector is  $\kappa(u) = (\rho - \lambda_5)^{-1}=5.5674\cdot 10^{-1}$.} 
The analysis of the matrix $R^{(\rho)}$ suggests equally that the edges 
$e(5 \leftrightarrow 4)$, $e(6 \leftrightarrow 4)$, or $e(6 \leftrightarrow 5)$ be
removed. This would weaken the community consisting of the nodes $\{4,5,6\}$. Edge
removal based on the Fiedler vector is not reliable, because the second smallest 
eigenvalue of the graph Laplacian matrix has multiplicity two: 
$\mu=\alpha_2=\alpha_3=0.6972$.

Consider instead the graph ${\cal G}$ that is shown in Figure \ref{FIG2} (right), in which
the edges that connect the three nodes that form { the community} have weight $0.5$. The
Perron value of ${\cal G}$ is $\rho=\lambda_{6}=1.6180$ and the Perron vector is 
\begin{equation*} 
u= \left[0.3035,\, 0.3035,\, 0.3035,\, 0.4991,\, 0.4991,\, 0.4991\right]^T.
\end{equation*}
{ One has $\kappa(u) = (\rho - \lambda_5)^{-1}=1.1944$.} 
The matrix $R^{(\rho)}$ suggests equally that the edges $e(4 \leftrightarrow 1)$, 
$e(5 \leftrightarrow 2)$, or $e(6 \leftrightarrow 3)$ be removed. This would weaken the 
anti-community. Hence, the situation is the opposite of the one above. Again, edge removal
based on the entries of the Fiedler vector is not possible since the Fiedler value is the
multiple eigenvalue $\mu=\alpha_2=\alpha_3=0.5$.

\begin{figure}
\centering
\begin{tabular}{cc}
\includegraphics[scale = 0.4]{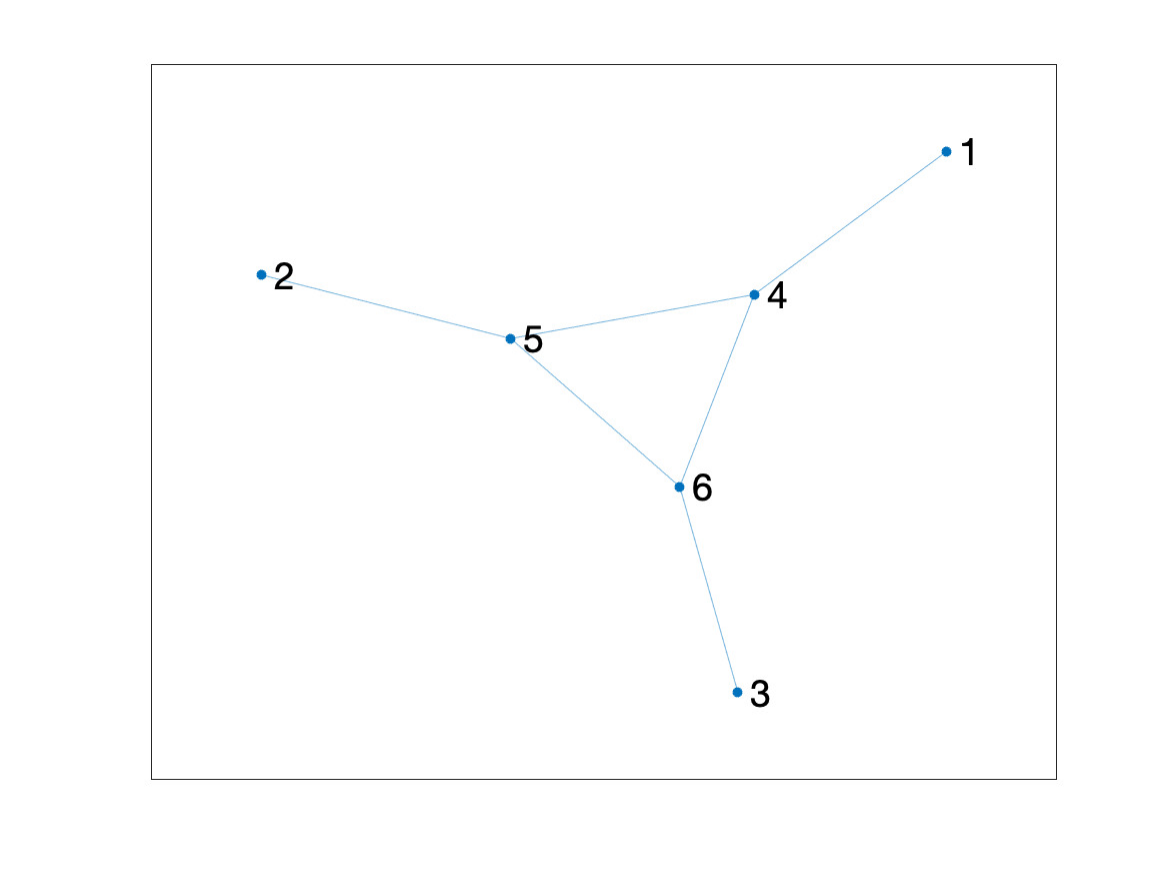} 
\includegraphics[scale = 0.4]{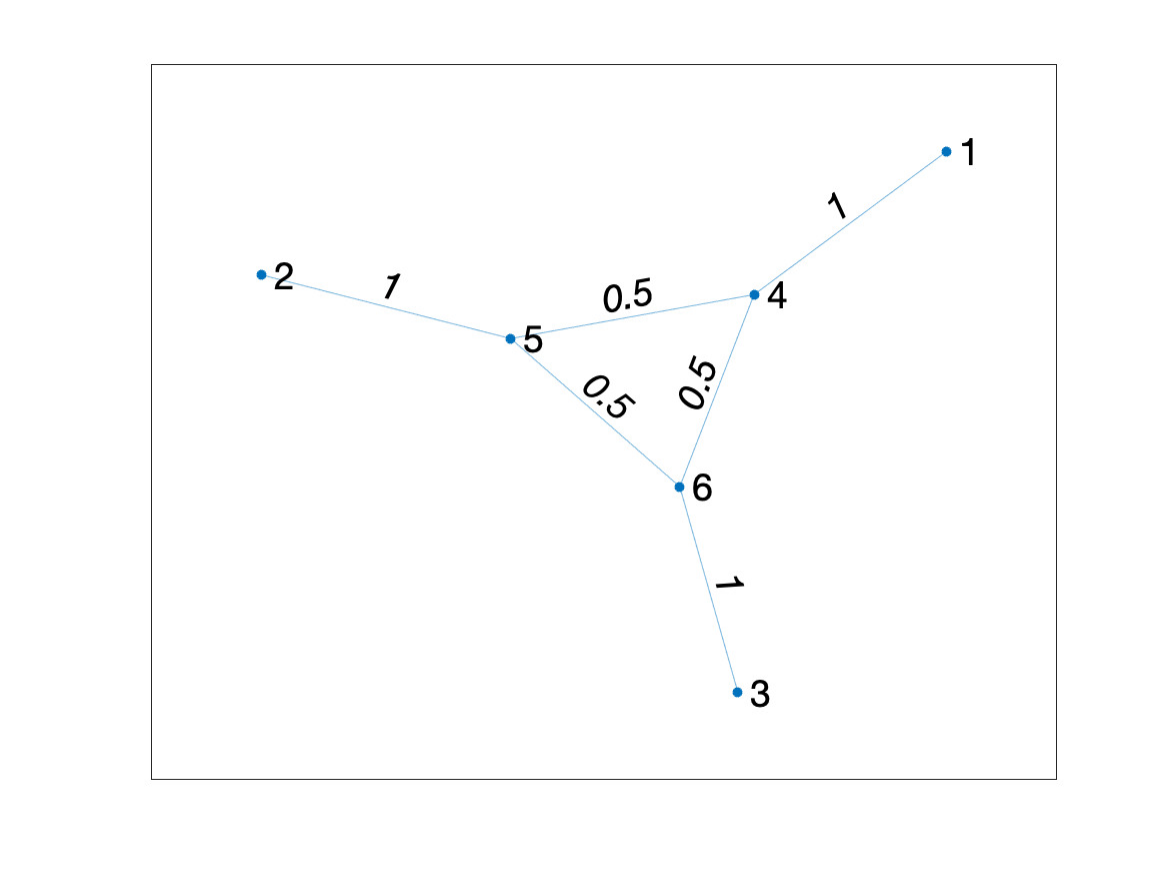} 
\end{tabular}
\caption{Unweighted graph  (left) and weighted graph (right) considered in Example \ref{ex3}.}
\label{FIG2}
\end{figure}

\end{example}

\begin{example}[Autobahn data set] \label{ex4}
This example considers an unweighted graph ${\cal G}$ that represents the German highway 
system network {\it Autobahn}. The graph is available at \cite{DCL}; see \cite{KH}. Its 1168 nodes are
German locations and its 1243 edges represent highway segments that connect these 
locations. 

The Perron value of the associated symmetric and irreducible adjacency matrix is 
$\rho=\lambda_{1168}=3.8251$ and the second smallest eigenvalue of the Laplacian matrix 
(the Fielder value) is $\mu=\alpha_2=0.0029$. The condition numbers of the Perron and 
Fiedler vectors are
\[
 \kappa(u)=(\rho-\lambda_{1167})^{-1}=6.7687 \mbox{~~~and~~~}
 \kappa(v)=(\alpha_3-\mu)^{-1}=1.2689\cdot 10^3,
\]
respectively. The edge whose removal has the largest impact on $\rho$ is 
$e(219 \leftrightarrow 217)$ and corresponds to the highway segment that connects Duisburg
(node $219$) and D\"usseldorf (node $217$). The edge whose removal has the largest impact 
on $\mu$ is $e(853 \leftrightarrow 850)$ and corresponds to the highway segment that 
connects Refrath (node $853$) and R\"udersdorf (node $850$). In detail, 
\begin{itemize}
\item[-] The entry $(219,217)$ of $R^{(\rho)}$ is $8.5961\cdot 10^{-2}$, while the average
of the positive entries of $R^{(\rho)}$ is $8.0580\cdot 10^{-4}$;
\item[-] The entry $(853,850)$ of $R^{(\mu)}$ is $4.7215\cdot 10^{-2}$, while the average 
of the positive entries of $R^{(\mu)}$ is $8.0451\cdot 10^{-4}$.
\end{itemize} 
Table \ref{T_2} shows the components of the Perron vector $u$ and the Fiedler vector $v$
associated with these nodes and their neighbors. The median of the Fiedler vector is 
$\beta=-7.3234\cdot 10^{-3}$. 

Let the graph ${\cal G}_1$ be obtained from $\cal G$ by removing edge 
$e(219 \leftrightarrow 217)$. Then ${\cal G}_1$ has the smaller Perron value $3.6819$ 
and (about the same) algebraic connectivity $0.0029$. In fact, after removal of this edge, 
Duisburg (node $219$) is connected with Dinslaken (node $198$), Dortmund (node $207$),  
Eiten (node $239$), Essen (node $267$), Flughafen (node $295$), and Krefeld (node $565$). 
These $6$ neighbors have on average $5$ neighbors. Moreover, D\"usseldorf (node $217$) 
remains connected with Diemelstadt (node $190$), Dinslaken (node $198$), Dorsten (node 
$206$), Dremmen (node $211$), Elmpt (node $236$), Enste (node $251$), and Erwitte 
(node $263$). These $7$ neighbors have on average $2.6$ neighbors. 


The graph ${\cal G}_2$, which is obtained from $\cal G$ by removing edge 
$e(853\leftrightarrow 850)$ has about the same Perron value  $3.8251$ and smaller 
algebraic connectivity $8.9349\cdot 10^{-4}$ as ${\cal G}_1$. In fact, after removal of
this edge, Refrath (node $853$) is only connected to Reichshoft (node $861$), which is 
connected only to Ronneburg (node $889$). As for R\"udersdorf (node $850$), it is 
connected to Pulsnitz (node $834$), Rangsdorf (node $842$), and Schwante (node $947$),
which are in turn also connected to only one other location each: Overath (node  $789$), 
Postdam (node $829$), and Velten (node $1048$), respectively.

We may observe { in Table \ref{T_2} that the nodes that represent  Duisburg, 
D\"usseldorf, and their neighbors correspond to (negative) components of the Fiedler vector
very close to the median $\beta$. On the contrary,} the nodes that 
represent Refrath, Reichshoft, their neighbors, and the 
neighbors of their neighbors correspond to positive components of the Fiedler vector that 
are an order of magnitude larger than $\beta$. 
Hence, these nodes most likely would be in 
the same subgraph defined by the Fiedler bipartition of $G$ (along with all their 
connections).

\begin{table}[h!tb]
\centering
\small
\begin{tabular}{c l c  c | c l c c }

$k$ & label           & $u_k$ & ${ v_k}$&$k$ & label           & $u_k$ & ${ v_k}$\\
\hline
$\bf{219}$  & Duisburg & $4.4247\cdot 10^{-1}$ & ${ -7.4852\cdot 10^{-3}}$ & 
$\bf{217}$  & D\"usseldorf  & $3.7156\cdot 10^{-1}$ & ${ -7.1928}\cdot 10^{-3}$ \\
$198$  & Dinslaken & $2.4965\cdot 10^{-1}$ & ${-7.4913\cdot 10^{-3}}$ & 
$190$  & Diemelstadt  & $1.0488\cdot 10^{-1}$ & ${ -6.9965}\cdot 10^{-3}$ \\
$207$  & Dortmund & $2.4746\cdot 10^{-1}$ & ${ -7.4741\cdot 10^{-3}}$ 
& $198$  & Dinslaken & $2.4965\cdot 10^{-1}$ & ${ -7.4913}\cdot 10^{-3}$ \\
$239$  & Eiten & $1.2519\cdot 10^{-1}$ & ${-7.7385\cdot 10^{-3}}$ & 
$206$  & Dorsten & $1.2623\cdot 10^{-1}$ & ${ -7.1365}\cdot 10^{-3}$ \\
$267$  & Essen & $3.2937\cdot 10^{-1}$ & ${ -7.4748\cdot 10^{-3}}$ & 
$211$  & Dremmen & $1.0488\cdot 10^{-1}$ & ${ -7.2992}\cdot 10^{-3}$ \\
$295$  & Flughafen & $1.6445\cdot 10^{-1}$ & ${-8.5000\cdot 10^{-3}}$ & 
$236$  & Elmpt & $1.8325\cdot 10^{-1}$ & ${ -7.3444}\cdot 10^{-3}$ \\
$565$  & Krefeld & $2.0483\cdot 10^{-1}$ & ${ -6.5037\cdot 10^{-3}}$ 
& $251$  & Enste & $1.0488\cdot 10^{-1}$ & $ { -6.9392}\cdot 10^{-3}$ \\
  & &  & &                                                            
  $263$  & Erwitte & $1.0502\cdot 10^{-1}$ & ${ -6.8299}\cdot 10^{-3}$ \\
\hline
$\bf{853}$  &  Refrath  & $5.1349\cdot 10^{-9}$ & $7.1819\cdot 10^{-2}$ 
& $\bf{850}$  & R\"udersdorf  & $1.8192\cdot 10^{-8}$ & $6.0193\cdot 10^{-2}$ \\
$861$  & Reichshoft & $1.4495\cdot 10^{-9}$ & $8.3240\cdot 10^{-2}$ 
& $834$  & Pulsnitz & $1.7770\cdot 10^{-8}$ & $5.3255\cdot 10^{-2}$ \\
  & &  & &                                                              
$842$  & Rangsdorf  & $4.1576\cdot 10^{-8}$ & $5.4986\cdot 10^{-2}$ \\
    & &  & &                                                            
$947$  & Schwante & $5.1048\cdot 10^{-9}$ & $6.0540\cdot 10^{-2}$ \\
\end{tabular}
\caption{Example \ref{ex4}. Spectral measures for selected nodes of the network 
{\it Autobahn}.} 
\label{T_2}
\end{table}
\begin{figure}
\centering
\begin{tabular}{cc}
\includegraphics[scale = 0.4]{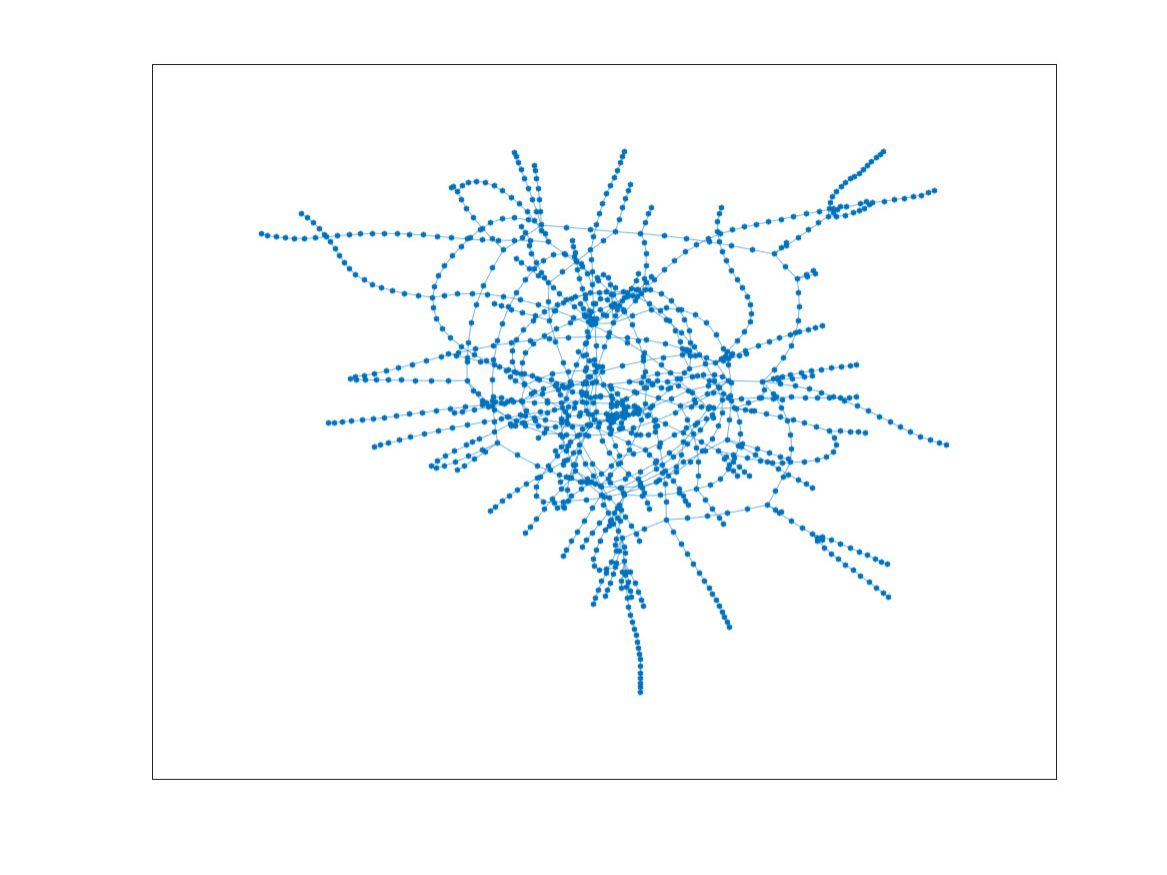} 
\includegraphics[scale = 0.4]{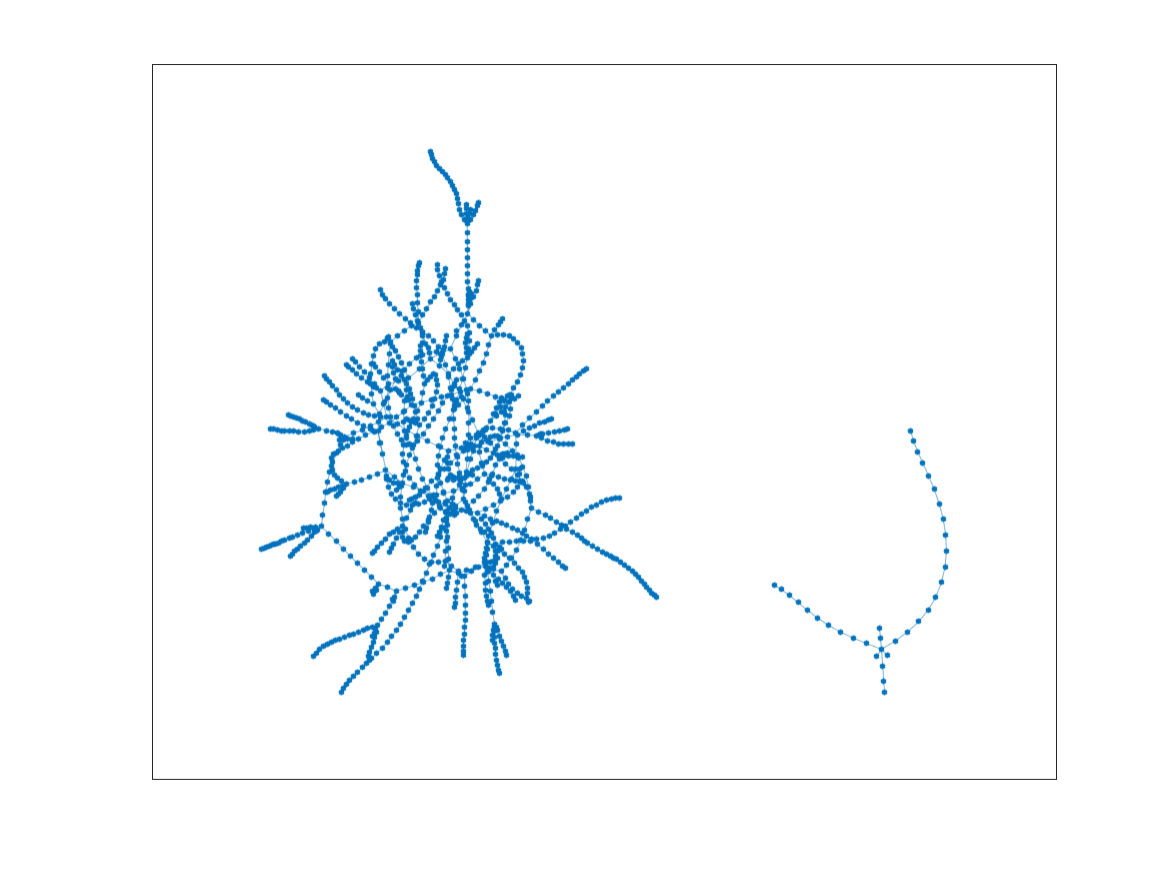} 
\end{tabular}
\caption{Example \ref{ex4}: Network {\it Autobahn} before (left) and after (right) 
removals suggested by the Perron vector.}
\label{FIG3}
\end{figure}

The bipartization of $\cal{G}$ achieved with the following procedure gives rise to the 
subgraphs in Figure \ref{FIG3} (right):
\begin{enumerate}
\item Find indices $(i,j)$ with $i>j$, such that $(i,j)=\arg\max R^{(\rho)}$ { for $A$};
\item Remove edge $e(i \leftrightarrow j)$ { by setting entries $a_{ij}$ and $a_{ji}$ to zero};
\item Goto step 1. and repeat until the graph is disconnected, { i.e., until $A$ becomes reducible}.
\end{enumerate}

Note that the components of the Fiedler vector $v$ of the graph $\cal{G}$ that are
associated to the nodes that after application of the above procedure live in the right 
subgraph of Figure \ref{FIG3} { (right)} all are negative. The spectral radius of the resulting 
adjacency matrix is $\hat\rho=\hat\lambda_{1168}={ 2.7926}$  and the { computed} second smallest 
eigenvalue of the resulting Laplacian matrix is $\hat\alpha_2={ 9.4687\cdot 10^{-17}}$.
{ The largest condition number of the computed Perron vectors is $8.4587 \cdot 10^{2}$.}
{    We conclude by observing that if the bipartization had been achieved without recomputing the matrix 
$R^{(\rho)}$ at each step, one would have obtained $\hat\rho=\hat\lambda_{1168}={ 3.6789}$  and $\hat\alpha_2={-1.7679\cdot 10^{-16}}$.}
\end{example}

\section{Some remarks on the multiplex case}\label{sec4}
Modeling complex systems that consist of $\ell$ different types of objects may lead to 
multilayer networks that can be represented by $\ell$ graphs that share the same set of 
$n$ nodes with nodes connected by edges both within a layer (intra-layer edges) and 
between layers (inter-layer edges).

We are concerned with multilayer networks that are \emph{multiplex} graphs, in which all 
the inter-layer edges have the same fixed weight $\gamma>0$. These edges determine the 
ease of communication between layers, and they connect each node with its copy in 
another {    layer}. A large value of $\gamma$ indicates strong interlayer coupling. 
Applications of multiplex graphs include modeling transportation networks; see, e.g., 
\cite{BS,NR2}.

A multiplex with $n$ nodes and $\ell$ layers may be represented by a supra-adjacency
matrix $B(\gamma)\in\mathbb{R}^{n\ell\times n\ell}$, which is a block matrix with 
$n\times n$ blocks of order $\ell$, see \cite{DSOGA},
\begin{equation*}
B=B(\gamma)={\rm blkdiag}[A^{(1)},A^{(2)},\dots,A^{(\ell)}]+\gamma
(1_\ell 1_\ell^T\otimes I_n-I_{n\ell}),
\end{equation*}
where the $h$th intra-layer adjacency matrix $A^{(h)}$ is the $h$th diagonal block, for 
$h=1,2,\ldots,\ell$, and each off-diagonal block in position $(h_1,h_2)$, with 
$1\leq h_1,h_2\leq \ell$ and $h_1\ne h_2$ equals $\gamma \, I_n$ and represents the
inter-layer connections between the layers $h_1$ and $h_2$. The operator $\otimes$ stands 
for the Kronecker product.

The supra-graph-Laplacian matrix $M(\gamma)\in\R^{n\ell\times n\ell}$ associated with the
multiplex, defined by
\[
M(\gamma)={\rm diag}[B(\gamma) 1_{n\ell}]-B (\gamma),
\]
is a block matrix with $n\times n$ blocks:
\begin{equation*}
M=M(\gamma)={\rm blkdiag}[L^{(1)},L^{(2)},\dots,L^{(\ell)}]+\gamma
(\ell I_{n\ell}-1_\ell 1_\ell^T\otimes I_n),
\end{equation*}
where the $h$th intra-layer graph Laplacian matrix 
$L^{(h)}= {\rm diag}[A^{(h)}1_{n}]-A^{(h)}$ is the $h$th diagonal block for 
$h=1,2,\ldots,\ell$.

Consider the Perron vector $u_B$ of $B(\gamma)$. In the context of eigenvector centrality 
metric, we refer to the matrix $U(\gamma)\in \R^{n\times \ell}$, formed column by column by
the entries of $u_B$, as the \emph{eigentensor} of $B$. The entry $(i,h)$ of $U(\gamma)$ is
said to be the \emph{joint eigenvector centrality} of node $i$ in layer $h$, and the $i$th 
entry of the vector $U(\gamma) 1_\ell$ is referred to the {\it eigenvector versatility} of
node $i$; see \cite{DSOGA}.
 
As for the Fiedler value $\mu_B$ of $B(\gamma)$, one has that if $\gamma\ll1$ (which
corresponds to weak interlayer coupling), then $\mu_B$ can be approximated by the smallest
positive eigenvalue of the inter-layer component of the supra-Laplacian matrix, 
$\gamma(\ell I_{n\ell}-1_{\ell}1_{\ell}^T\otimes I_{n})$, that is by $\gamma \ell$. On the
contrary, if $\gamma\gg1$, then $\mu_B$ can be approximated by the second eigenvalue 
$\mu_{\bar A}$ of the average network $\bar A=\frac{1}{\ell}\sum_{h=1}^{\ell} A^{(h)}$; 
see \cite{SDKDGA}. In general, 
\[
\mu_B\leq \min\{\mu_{\bar A}, \gamma \ell\}.
\]

\subsection{A numerical example}
The considerations regarding the single-layer Perron vector procedure are difficult to 
extend to the multiplex case, while the single-layer Fiedler procedure cannot be 
extended to the multiplex case at all. What can be done is to treat aggregate adjacency 
matrices associated with multiplexes. We illustrate this in the following example. 

\begin{example}[The Tube data set] \label{ex6}
In this example we consider the (weighted) adjacency matrix $W=[w_{ij}]\in\R^{n\times n}$ 
that is associated with the $n=315$ stations of the London Underground Transportation 
Network, {\it The Tube}. The weight $w_{ij}=w_{ji}$ is set to be the number of contiguous 
connections between stations $i$ and $j$. Thus, $W$ is an aggregated adjacency 
matrix associated with the (unweighted) multiplex formed by $\ell=14$ layers (subway 
lines Bakerloo, Central, Circle, District, Hammersmith $\&$ City, Jubilee, Metropolitan, 
Northern, Piccadilly, Victoria, Waterloo, Overground, DLR, and DLR New). Thus, one has 
$W=\sum_{h=1}^{\ell} A^{(h)}$, where the $h$th intra-layer adjacency matrix $A^{(h)}$ is 
for line $h$. The Perron value of the network is $\rho= 6.1827$ and the Perron vector has
condition number $\kappa(u)=1.4500$.

The edge $e(148 \leftrightarrow 92)$ with weight $3$ has the largest impact on $\rho$.  
There are $3$ lines that connect King's Cross St. Pancras (node $148$) and Farringdon 
(node $92$), i.e., Circle (layer $3$), Hammersmith $\&$ City (layer $5$), and Metropolitan
(layer $7$). Decreasing the weight of the edge $e(148\leftrightarrow 92)$ by one, the 
Perron value decreases to $5.8801$.

The following procedure yields the bipartition of the network depicted in Figure 
\ref{FIG4} (left) shown in Figure \ref{FIG4} (right) with Perron value $2.3468$.

\begin{enumerate}
\item Find indices $(i,j)$ with $i>j$ such that $(i,j)=\arg\max R^{(\rho)}$ for $W$;
\item Decrease both the entries $w_{ij}$ and $w_{ji}$ by $1$;
\item Goto step 1. and repeat until the  graph is disconnected.
\end{enumerate}

\begin{figure}
\centering
\begin{tabular}{cc}
\includegraphics[scale = 0.4]{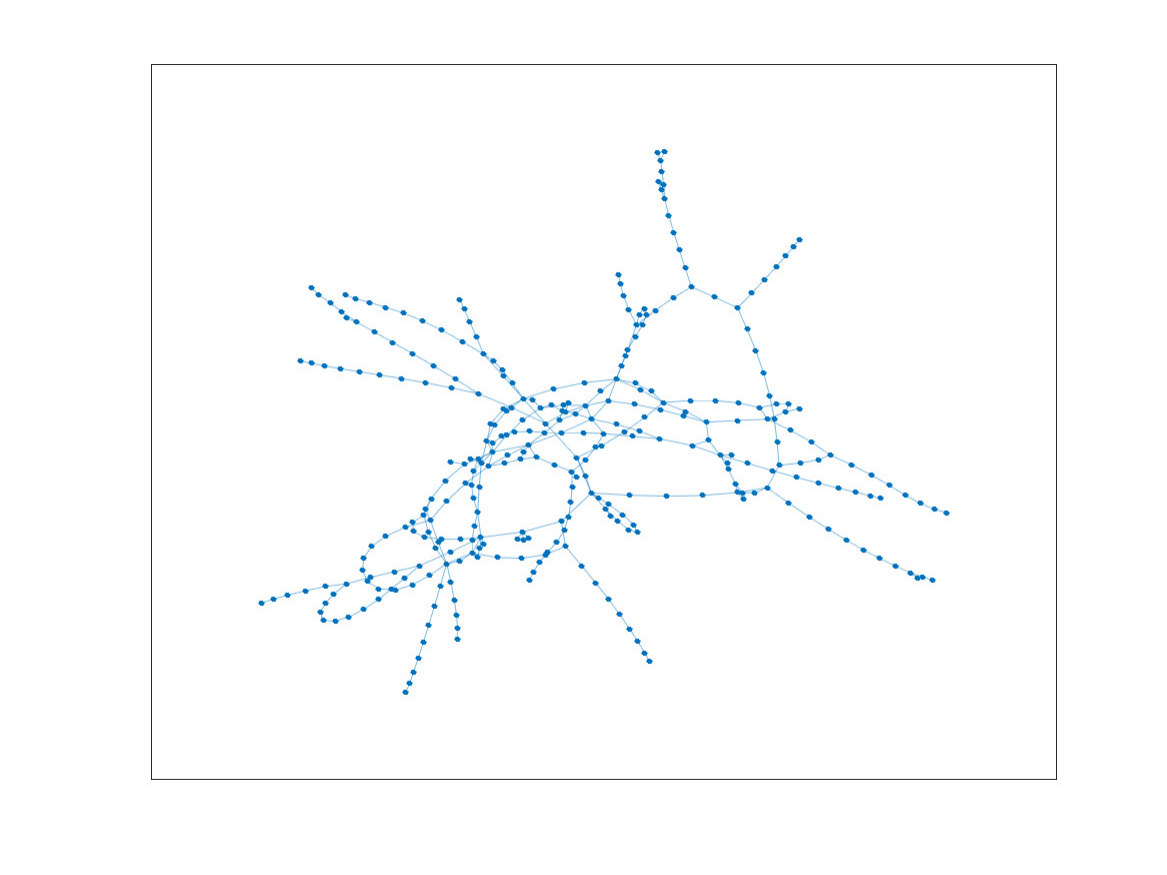} 
\includegraphics[scale = 0.4]{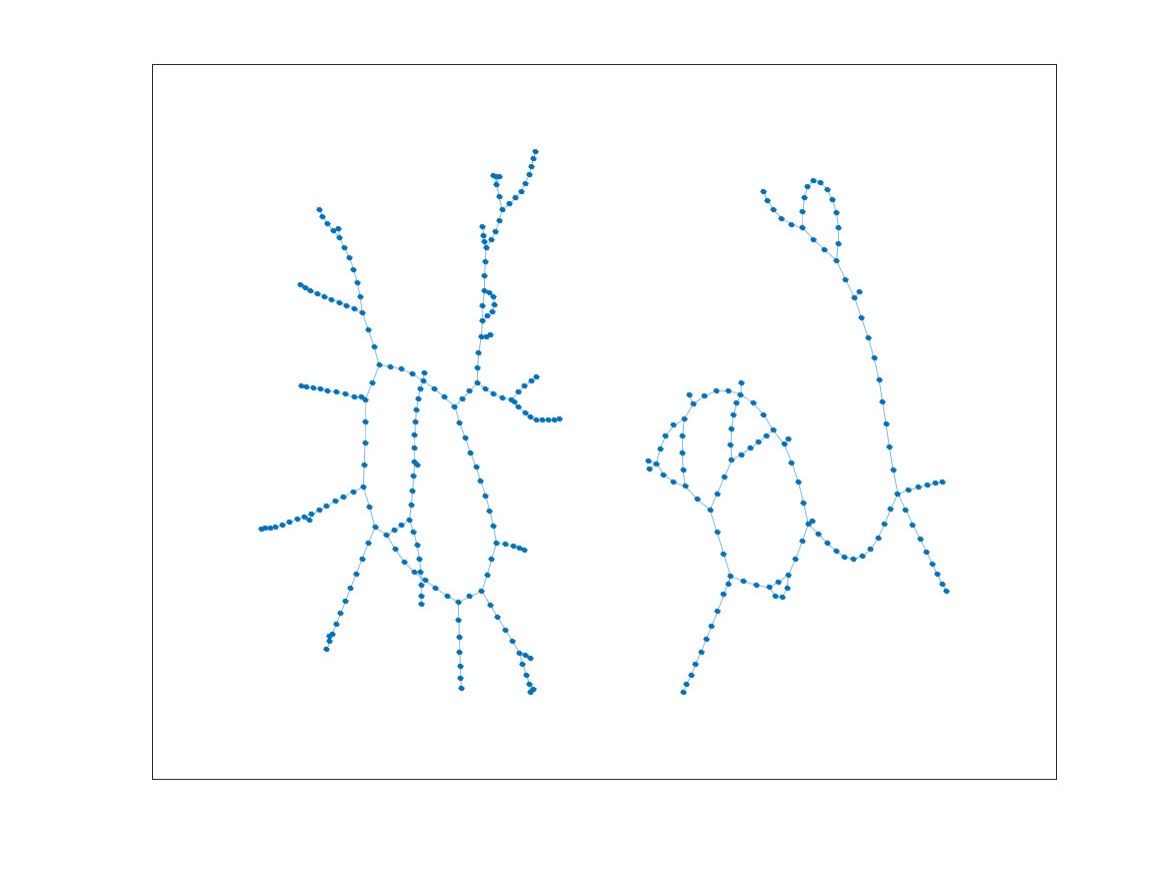} 
\end{tabular}
\caption{Example \ref{ex6}: Aggregated network {\it The Tube}, before (left) and after 
(right) edge removals suggested by the Perron vector.}
\label{FIG4}
\end{figure}
\end{example}

\section{Conclusions}\label{sec5}
In this paper we investigated the sensitivity of Perron and Fiedler values to edge-weight 
perturbations as well as the conditioning of their corresponding eigenvectors. This 
comprehensive perturbation analysis is applied to the identification of edges that
are critical to the 
structural robustness of the network. We observed that the Perron sensitivity analysis 
tends to select central edges, whose removal would facilitate beneficial network 
simplification and give structurally robust sub-networks. 
Our findings illustrate that Perron and Fiedler sensitivity analyses offer complementary 
insights into network structure and robustness. 

\section*{Acknowledgments}
{    The authors would like to thank the referees for comments.}
Research by SN was partially supported by a grant from SAPIENZA Universit\`a di Roma, 
by INdAM-GNCS, and by the European Union – NextGenerationEU CUP F53D23002700006: PRIN 2022 
research project “Inverse Problems in the Imaging Sciences (IPIS)", grant n. 2022ANC8HL. 


\begin{thebibliography}{99}
\bibitem{ABCMP} 
{    D. Altafini, D. A. Bini, V. Cutini, B. Meini, and F. Poloni, An edge centrality 
measure based on the Kemeny constant, SIAM Journal on Matrix Analysis and Applications, 
44 (2023), pp. 648--669.}
\bibitem{AB} 
{    F. Arrigo and M. Benzi, Updating and downdating techniques for optimizing
network communicability, SIAM Journal on Scientific Computing, 38 (2016), pp. B25--B49.}
\bibitem{BCR1}
J. Baglama, D. Calvetti, and L. Reichel, IRBL: An implicitly restarted block Lanczos 
method for large-scale Hermitian eigenproblems, SIAM Journal on Scientific Computing, 24 
(2003), pp. 1650--1677.
\bibitem{BCR2}
J. Baglama, D. Calvetti, and L. Reichel, Algorithm 827: irbleigs: a MATLAB program for 
computing a few eigenpairs of a large sparse Hermitian matrix, ACM Transactions on 
Mathematical Software, 29 (2003), pp. 337--348.
\bibitem{BS}
K. Bergermann and M. Stoll, Orientations and matrix function-based centralities in 
multiplex network analysis of urban public transport, Applied Network Science, 6 (2021), 
pp. 1--33.
\bibitem{PF}
A. Berman and R. J. Plemmons, Nonnegative Matrices in the Mathematical Sciences, SIAM, 
Philadelphia, 1994.  
\bibitem{Bi}
N. Biggs, Algebraic Graph Theory, Cambridge University Press, Cambridge, 1994. 
\bibitem{Bo} 
P. Bonacich, Power and centrality: a family of measures, American Journal of Sociology, 92
(1987), pp. 1170--1182.
\bibitem {Br}
{    D. Braess, \"Uber ein Paradoxon aus der Verkehrsplanung, Unternehmensforschung, 12 
(1968), pp. 258--268.}
\bibitem{BH} 
A. E. Brouwer and W. H. Haemers, Spectra of Graphs. Springer, Heidelberg, 2012. 
\bibitem{BC} 
S. Butler and F. Chung, Spectral Graph Theory, in Handbook of Linear Algebra, 2nd edition,
Leslie Hogben, ed., CRC Press, 2014.
\bibitem{Ch} 
F. Chung, Spectral graph theory, CBMS Regional Conference Series in Mathematics, vol. 92,
American Mathematical Society, Providence, 1997.
\bibitem{CNRR}
A. Concas, S. Noschese, L. Reichel, and G. Rodriguez, A spectral method for bipartizing a 
network and detecting a large anti-community, Journal of Computational and Applied 
Mathematics, 373 (2020), Art. 112306.  
\bibitem{DCL}
Dynamic Connectome Lab - Data Sets. \hfill\break
\url{https://sites.google.com/view/dynamicconnectomelab}
\bibitem{DSOGA}
M. De Domenico, A. Sol\'e-Ribalta, E. Omodei, S. G\'omez, and A. Arenas, Centrality in 
interconnected multilayer networks, arXiv:1311.2906v1  (2013).
\bibitem{DCMR}
{    O. De la Cruz Cabrera, M. Matar, and L. Reichel, Centrality measures for 
node-weighted networks via line graphs and the matrix exponential, Numerical Algorithms, 
88 (2021), pp.  583--614.}
\bibitem{Es}
E. Estrada, The Structure of Complex Networks: Theory and Applications, Oxford
University Press, Oxford, 2011.
\bibitem{FT}
D. Fasino and F. Tudisco, A modularity based spectral method for simultaneous community 
and anti-community detection, Linear Algebra and Its Applications, 542 (2017), pp. 
605--623. 
\bibitem{F1}
M. Fiedler, Algebraic connectivity of graphs, Czechoslovak Mathematical Journal, 23 (1973),
no. 98, pp. 298--305. 
\bibitem{F2}
M. Fiedler,  A property of eigenvectors of nonnegative symmetric matrices and its 
application to graph theory, Czechoslovak Mathematical Journal, 25 (1975), no. 4, 
pp. 619--633. 
\bibitem{GBS}
{    A. Ghosh, S. P. Boyd, and A. Saberi, Minimizing effective resistance of a graph, SIAM
Review, 50 (2008), pp. 37--66.}
\bibitem{HJ}
R. A. Horn and C. R. Johnson, Matrix Analysis, Cambridge University Press, Cambridge, 
1985. 
\bibitem{HJ2}
R. A. Horn and C. R. Johnson, Topics in Matrix Analysis, Cambridge University Press, 
Cambridge, 1991. 
\bibitem{JKVV} 
A. Jamakovic, R. E. Kooij, P. Van Mieghem, E. R. van Dam, Robustness of networks against 
viruses: the role of the spectral radius, in Symposium on Communications and Vehicular 
Technology, 2006, pp. 35--38.
 \bibitem{KH} 
{\ M.  Kaiser and C. C. Hilgetag, Spatial growth of real-world networks, Physical 
Review E, 69 (2004), Art. 036103.}
 \bibitem{KKT}
{ M. Karow, D. Kressner, and F. Tisseur, Structured eigenvalue condition numbers, SIAM 
Journal on Matrix Analysis and Applications, 28 (2006), pp. 1052--1068.}
\bibitem{MT}
{    S. Massei and F. Tudisco, Optimizing network robustness via Krylov subspaces, 
ESAIM Mathematical Modelling and Numerical Analysis. 58 (2024), pp. 131--155.}
\bibitem{VM}
P. V. Mieghem, Graph Spectra for Complex Networks, Cambridge University Press, Cambridge,
2011. 
\bibitem{MSN}
A. Milanese, J. Sun, T. Nishikawa, Approximating spectral impact of structural 
perturbations in large networks, Physical Review E, 81 (2010), Art. 046112.
\bibitem{Ne}
M. E. J. Newman, Networks: An Introduction, Oxford University Press, Oxford, 2010.
\bibitem{NP}
S. Noschese and L. Pasquini, Eigenvalue condition numbers: zero-structured versus 
traditional, Journal of Computational and Applied Mathematics, 185 (2006), pp. 174--189.
\bibitem{NR}
S. Noschese and  L. Reichel, Eigenvector sensitivity under general and structured 
perturbations of tridiagonal Toeplitz-type matrices, Numerical Linear Algebra with 
Applications, 26 (2019), no. 3, Art. e2232.
\bibitem{NR1}
S. Noschese and L. Reichel, Estimating and increasing the structural robustness of a 
network, Numerical Linear Algebra with Applications 29 (2021), no. 2, Art. e2418.
\bibitem{NR3}
{    S. Noschese and L. Reichel, Enhancing multiplex global efficiency, Numerical 
Algorithms, 96 (2024), pp. 397--416.}
\bibitem{NR4}
S. Noschese and L. Reichel, Edge importance in complex networks, Numerical Algorithms, 99 
(2025), pp. 377--410.
\bibitem{NR2}
S. Noschese and L. Reichel, Communication in multiplex transportation networks, Numerical 
Algorithms, in press. 
\url{DOI: 10.1007/s11075-024-01943-4}
\bibitem{Sc}
{    M. Schweitzer, Sensitivity of matrix function based network communicability
measures: computational methods and a priori bounds, SIAM Journal on Matrix Analysis and 
Applications, 44 (2023), pp. 1321--1348.}
\bibitem{SDKDGA}
A. Sol\'e-Ribalta, M. De Domenico, N. E. Kouvaris, A. D\'iaz-Guilera, S. Gómez, and 
A. Arenas, Spectral Properties of the Laplacian of multiplex networks, Physical Review E 
88  (2013), Art. 032807.
\bibitem{St01}
G. W. Stewart, Matrix Algorithms, vol. II: Eigensystems, SIAM, Philadelphia, 2001. 
\bibitem{Wi}
J. H. Wilkinson, Sensitivity of eigenvalues II, Utilitas Mathematica, 30 (1986), pp. 
243--286.
\end{thebibliography}
\end{document}